\documentclass[12pt]{amsart}
\usepackage{amssymb}
\usepackage[dvips]{graphics}
\ifx\mathrm\undefined\let\mathrm\rm\fi
\ifx\mathbf\undefined\let\mathbf\bf\fi
\ifx\mathfrak\undefined\let\mathfrak\frak\fi
\ifx\mathcal\undefined\let\mathcal\cal\fi
\ifx\mathbb\undefined\let\mathbb\Bbb\fi
\ifx\emph\undefined\let\emph\it\fi
 at9.98pt

\newcommand{\g}{{{\mathfrak g}\,}}

\newcommand{\n}{{{\mathfrak n}}}
\newcommand{\B}{{{\mathfrak b}}}
\newcommand{\h}{{{\mathfrak h\,}}}

\newcommand{\Z}{{\mathbb Z}}

\newcommand{\C}{{\mathbb C}}

\newcommand{\Ref}[1]{{(\ref{#1})}}

\newcommand{\la}{\lambda}

\newcommand{\dontprint}[1]
{\relax}

\newtheorem%
{thm}{Theorem}[section]
\newtheorem%
{proposition}[thm]{Proposition}
\newtheorem%
{lemma}[thm]{Lemma}
\newtheorem%
{lemmadef}[thm]{Lemma-Definition}
\newtheorem%
{corollary}[thm]{Corollary}
\newtheorem%
{conjecture}[thm]{Conjecture}

\newcommand{\bea}{\begin{eqnarray*}}
\newcommand{\eea}{\end{eqnarray*}}
\newcommand{\bean}{\begin{eqnarray}}
\newcommand{\eean}{\end{eqnarray}}

\newcommand{\nc}{\newcommand}
\nc{\on}{\operatorname}
\nc{\al}{\alpha}
\nc{\ri}{\rangle}
\nc{\lef}{\langle}
\nc{\W}{{\mathcal W}}
\nc{\La}{\Lambda}
\nc{\ep}{\epsilon}
\nc{\Om}{\Omega}
\newcommand{\be}{\begin{displaymath}}
\newcommand{\ee}{\end{displaymath}}
\newcommand{\bs}{\boldsymbol}
\nc{\PCr}{{ \bs P  (\C[x])^r   }}
\newtheorem{theorem}{Theorem}[section]
\newtheorem{conj}[theorem]{Conjecture}

\newcommand{\p}{{\partial}}
\newcommand{\tG}{{{}^tG}}
\newcommand{\tg}{{{}^t\g}}
\newcommand{\bP}{{\mathbb P^1}}
\newcommand{\tB}{{{}^tB_-}}
\newcommand{\id}{{\rm{id}}}
\newcommand{\om}{{\rm{Om}_D}}

\newcommand{\tN}{{{}^tN_+}}

\begin{document}
\title[Miura Opers and Critical Points]
{Miura Opers and Critical Points of Master Functions}

\author[{}]{ Evgeny Mukhin  
${}^{*,1}$ 
\and Alexander Varchenko${}^{**,2}$}
\thanks{${}^1$ Supported in part by NSF grant DMS-0140460}

\thanks{${}^2$ Supported in part by NSF grant DMS-0244579}

\begin{abstract}
Critical points of a master function associated to a simple Lie  algebra $\g$
come in families called the populations \cite{MV1}. We prove
that a population is isomorphic to the flag variety of the Langlands dual Lie algebra
${}^t\g$. The proof is based on the correspondence between critical points and differential
operators called the Miura opers.

For  a Miura oper $D$,  associated with a critical point of a
population, we show that all solutions of the differential equation $DY=0$ 
can be written explicitly in terms of critical points composing the population.

\end{abstract}

\maketitle 
\medskip \centerline{\it ${}^*$
Department of Mathematical Sciences,} 
\centerline{\it Indiana University 
Purdue University Indianapolis,}
\centerline{\it 402 North Blackford St., Indianapolis,
IN 46202-3216, USA}
% \newline mukhin@math.iupui.edu}
 \medskip
\centerline{\it ${}^{**}$Department of Mathematics, University of
  North Carolina at Chapel Hill,} \centerline{\it Chapel Hill, NC
  27599-3250, USA} \medskip

%\maketitle

%\centerline{December, 2003}

\section{Introduction}
In \cite{ScV} rational functions were considered which are
products of powers of linear functions. 
It was discovered that under certain conditions 
all critical points of the rational functions are non-isolated 
and form non-trivial varieties. It is
not clear yet how general that phenomenon is but the phenomenon
certainly holds for products of powers of linear functions appearing in representation
theory. Those products are called the master functions.

Let $\h$ be the Cartan subalgebra of a simple Lie algebra $\g$; $( \ ,
\ )$ the Killing form on $\h^*$; $\al_1, \dots , \al_r \in \h^*$
simple roots; $\La_1, \dots , \La_n \in \h^*$ dominant integral
weights; $l_1, \dots , l_r$ non-negative integers; $z_1, \dots , z_n$
distinct complex numbers. The master function $\Phi$ associated with
this data is given by formula \Ref{master}. It is a
 rational function of $l_1 + \dots + l_r$ variables $
t_1^{(1)}, \dots , t_{l_1}^{(1)}, t_1^{(2)}, \dots , t_{l_r}^{(r)}$, and
$n$ variables $z_1, \dots , z_n$.  We may think that $l_1 + \dots + l_r + n$
particles are given in $\C$. The particle $t^{(i)}_j$ has weight $-
\al_i$ and the particle $z_s$ has weight $\La_s$. The particles
interact pairwise.  The interaction of particles $x$ and $y$ with
weights $v$ and $w$, respectively, is given by $(x - y)^{(v,w)}$. Then
total interaction is the product of such terms over the set of all pairs. The
master function describes the interaction of $t$-particles 
 and $z$-particles.

The master function appears in hypergeometric solutions to the KZ
equations with values in the tensor product of irreducible highest
weight representations $L_{\La_1}, \dots , L_{\La_n}$ with highest
weights $\La_1, \dots , \La_n$, respectively.  The solutions have the
form 
\be 
u (z)\ = \ \int \ \Phi (\bs t; \bs z )^{1/ \kappa}\ A(
\bs t; \bs z ) \ d \bs t\ , 
\ee 
where $\kappa$ is the parameter of the
KZ equations and $A( \bs t; \bs z )$ is some explicitly written
rational function with values in the tensor product \cite{SV}.

The master function also appears in the Bethe ansatz of the Gaudin
model with values in the same tensor product \cite{RV}.  In that case
the value of the function $A(\ \cdot\ ; \bs z )$ at a point $\bs t$ is
an eigenvector of the commuting Gaudin Hamiltonians if $\bs t$ is a
critical point of the master function.

In this paper we study critical points of the master function on the
set where all $\{ t^{(i)}_j, z_s\}$ are distinct.  In other words we
study those positions of distinct particles $\{ t^{(i)}_j \}$ in the
complement to $\{z_s\}$ which extremize the master function.

Critical points of master functions associated to a simple Lie
algebra $\g$ come in families called populations \cite{ScV, MV1}. In this
paper we prove that a population is isomorphic to the flag variety of
the Langlands dual Lie algebra ${}^t\g$. The proof is based on the
correspondence between critical points and differential operators
called the Miura opers.

To every critical point $\bs t$ one assigns a certain linear
differential operator $D_{\bs t}$ with coefficients in ${}^t\g$,
called the Miura oper. The differential operators of that type were
considered by V. Drinfeld and V. Sokolov in their study of the KdV type
equations \cite{DS}. On opers and Miura opers see \cite{BD, FFR, F1, F2, MV1, 
BM, MV4}.

Different critical points correspond to different Miura opers. The
Miura opers corresponding to critical points of a given population
form an equivalence class with respect to suitable gauge
equivalence. We show that the equivalence class of Miura opers is
isomorphic to the flag variety of ${}^t\g$.

In \cite{MV1, BM} we considered Miura opers for Lie algebras of types 
$A_r, B_r, C_r, G_2$ and using 
the opers proved that a population of critical points of types 
$A_r, B_r, C_r, G_2$ is isomorphic to the 
corresponding flag variety. The proof, suggested in the present
 paper, is more direct and works for any simple Lie algebra. 

If $D_{\bs t}$ is the Miura oper corresponding to a critical point
$\bs t$, then the set of solutions of the differential equation
$D_{\bs t}\,Y\,=\,0$ with values in a suitable space is an important
characteristics of the critical point. We used that characteristics
for $sl_{r+1}$ Miura opers in \cite{MV1} to give a bound from above
for the number of populations of critical points  of the corresponding
$sl_{r+1}$ master function. That statement in \cite{MV1} was
in some sense opposite 
to the Bethe ansatz conjectures, see \cite{MV1}.

 It turns out that for any simple Lie algebra
$\g$ and any critical point $\bs t$ of a $\g$ master function all
solutions of the differential equation $D_{\bs t}\,Y\,=\,0$ can be
written explicitly in terms of critical points composing the
population originated at $\bs t$.  Thus the population of critical
points ``solves'' the Miura differential equation $D_{\bs t}\,Y\,
=\,0$. This is the second main result of the paper.

\bigskip

When this paper was being written preprint \cite{F2} by E. Frenkel
appeared. The preprint is devoted to the same fact that the variety of
gauge equivalent ${}^t\g$ Miura opers is isomorphic to the flag variety of
${}^t\g$. The main result of \cite{F2} is Corollary 3.3. In our opinion the
proof of Corollary 3.3 in \cite{F2} is incomplete and in particular the 
statement of Corollary 3.3 is not correct.

The idea of this paper was originated in discussions with E. Frenkel
in the spring of 2002. As a result of those discussions two
papers appeared: this one and \cite{F2}.

We thank E. Frenkel for stimulating meetings which originated this paper.  
We thank P. Belkale and S. Kumar for numerous useful discussions.

\bigskip

The paper is organized as follows. In Section 2 we introduce populations 
of critical points. In Section 3 we discuss elementary properties of Miura opers
corresponding to critical points. In Section 4 we prove that the variety of 
gauge equivalent Miura opers is isomorphic to the flag variety, see Theorem
\ref{isomorphism}. We discuss the relations between the Bruhat cell decomposition of the
flag variety and  populations of critical points in Section 5. The main result there is
Corollary \ref{main corollary} describing the structure of connected 
components of the critical set of master functions.
In Section 6 we give explicit formulas for solutions of the differential
equation $D_{\bs t}\,Y\,=\,0$, see Theorems \ref{sl}, \ref{bc}, \ref{Solutions}.

\section{Master functions and critical points, \cite{MV1}}\label{crit pts}
\label{crit point sec}
\subsection{Kac-Moody algebras}\label{Kac_Moody sec}
Let $A=(a_{i,j})_{i,j=1}^r$ be a generalized  Cartan matrix, 
$a_{i,i}=2$,
$a_{i,j}=0$ if and only $a_{j,i}=0$,
 $a_{i,j}\in \Z_{\leq 0}$ if $i\ne j$. 
We  assume that $A$ is symmetrizable, i.e. 
there exists a diagonal matrix $D=\on{diag}\{d_1,\dots,d_r\}$ 
with positive integers $d_i$ such that $B=DA$
is symmetric.

Let $\g=\g(A)$ be the corresponding complex Kac-Moody 
Lie algebra (see \cite{K}, \S 1.2), 
$\h \subset \g$  the Cartan subalgebra.
The associated scalar product is non-degenerate on $\h^*$ and 
 $\on{dim}\h = r + 2d$,  where $d$ is the dimension of 
the kernel of the Cartan matrix $A$.

Let $\al_i\in \h^*$, $\al_i^\vee\in \h$, $i = 1, \dots , r$, be the sets of simple roots,
coroots, respectively. We have
\bea
 (\al_i,\al_j)&=& d_i \ a_{i,j}, \\
\langle\la ,\al^\vee_i\rangle&=&2(\la,\al_i)/{(\al_i,\al_i)}, \qquad \la\in\h^*.
\eea
In particular, $\langle\al_j ,\al^\vee_i \rangle = a_{i,j}$.

Let $\mathcal P = \{ \lambda \in \h^* \, |\, \langle\la
,\al^\vee_i\rangle \in \Z\}$ and $\mathcal P^+ = \{ \lambda \in \h^* \, |\, \langle\la
,\al^\vee_i\rangle \in \Z_{\geq 0}\}$ be the sets of integral and
dominant integral weights.

Fix $\rho\in\h^*$ such that $\langle\rho,\al_i^\vee\rangle=1$,
$i=1,\dots,r$. We have $(\rho,\al_i)= (\al_i,\al_i)/2$.

The Weyl group $\mathcal W\in\on{End (\h^*)}$ is generated by 
reflections $s_i$, $i=1,\dots,r$, 
\be
s_i(\la)=\la-\langle\la,\al_i^\vee\rangle\al_i, \qquad \la\in\h^*.
\ee
We use the notation
\bea\label{shifted}
w\cdot\la=w(\la+\rho)-\rho,\qquad w\in \mathcal W,\;\la\in\h^*,
\eea
for the shifted action of the Weyl group.

The Kac-Moody algebra $\g(A)$ is generated by $\h$, $e_1, \dots , e_r, f_1,
\dots , f_r$ with defining relations
\bea
[e_i, f_j] & = & \delta_{i,j} \,\alpha_i^\vee ,
\qquad i, j = 1, \dots r ,
\\
{}[ h , h'] & = & 0 ,
\qquad h, h' \in \h ,
\\
{}[ h, e_i] &=& \langle \alpha_i, h \rangle\, e_i ,
\qquad  h \in \h, \ i = 1, \dots r ,
\\
{}[ h, f_i] &=& - \langle \alpha_i, h \rangle\, f_i ,
\qquad  h\in \h, \ i = 1, \dots r ,
\eea
and the Serre's relations
\bea
(\mathrm{ad}\,{} e_i)^{1-a_{i,j}}\,e_j = 0 ,
\qquad
(\mathrm{ad}\, {} f_i)^{1-a_{i,j}}\,f_j = 0 ,
\eea
for all $i\neq j$. The generators
$\h$, $e_1, \dots , e_r, f_1,
\dots , f_r$ are called the Chevalley generators.

Denote  $\n_+$ (resp. $\n_-$) the subalgebra
generated by $e_1, \dots , e_r$ (resp. $f_1, \dots , f_r$). Then
$\g = \n_+\oplus \h\oplus \n_-$. Set $\B_\pm = \h \oplus \n_\pm$.

Let $\g = \oplus_j \g^j$ be the canonical grading of $\g$. Here we have
$e_i \in \g^1,\, f_i \in \g^{-1},\, \n_+ = \oplus_{j>0}\, \g^j,\,
\h = \g^0,\,\n_- = \oplus_{j<0}\, \g^j$.
%Let $\n_+^\vee$ be the completion of $\n_+$ with respect to the canonical grading.

For a vector space $X$ we denote $M(X)$ the space of $X$-valued rational functions on $\C$.

We denote  $\bar M(\n_+)$  the completion of the space $M(\n_+)$ with respect to 
the canonical grading. An element of $\bar M(\n_+)$ is 
a formal sum $\sum_{j>0} u_j$,
where $u_j : \C \to \g^j$ are rational functions.

The Kac-Moody algebra $^t\g=\g( ^tA)$ 
corresponding to the transposed Cartan matrix $^tA$ is called
{\it Langlands dual} to $\g$. 
Let $^t\al_i\in {}^t\h^*$, $^t\al_i^\vee\in {}^t\h$, 
$i = 1, \dots , r$, be the sets of simple roots,
coroots of $^t\g$, respectively. Then 
\bea
\langle {}^t\al_i, {}^t\al_j^\vee\rangle  =  \langle\al_j,\al_i^\vee\rangle = a_{i,j}
\eea
 for all $i,j$.

\subsection{The definition of
master functions and critical points}\label{master sec} 
We fix a Kac-Moody algebra
$\g=\g(A)$, a non-negative integer $n$, a collection of dominant integral weights  
$\bs\La=(\La_1, \dots , \La_n)$, $\La_i\in\mathcal P^+$,
and points $\bs z = \{z_1, \dots , z_n\} \subset \C$. We assume that
$z_i\neq z_j$ if $i\neq j$. We often do not stress the dependence of our
objects on these parameters.

In addition we choose 
a collection of non-negative integers  $\bs l=(l_1,\dots,l_r)\in\Z^r_{\geq 0}$. 
The choice of $\bs l$ is equivalent to the choice of the weight 
\be
\La_\infty\ =\sum_{i=1}^n \La_i  -  \sum_{j=1}^r
l_j\al_j \in \mathcal P .
\ee
The weight $\La_\infty$ will be called {\it the weight at infinity}.

The {\it master function} $\Phi(\bs t;\La_\infty)$
is defined by 
\bean\label{master}
&&
\Phi(\bs t; \La_\infty) = \Phi(\bs t; \bs z, \bs \La, \La_\infty) = 
\prod_{1\leq s < u \leq n}(z_s - z_u)^{(\La_s,\La_u)} \times
\\
&& 
\prod_{i=1}^r\prod_{j=1}^{l_i}\prod_{s=1}^n
  (t_j^{(i)}-z_s)^{-(\La_s,\al_i)}
  \prod_{i=1}^r\prod_{1\leq j<s\leq l_i} (t_j^{(i)}-t_s^{(i)})^{(\al_i,\al_i)}
  \prod_{1\leq i<j\leq r}
\prod_{s=1}^{l_i}\prod_{k=1}^{l_{j}}(t_s^{(i)}-t_k^{(j)})^{(\al_i,\al_j)},
\notag
\eean
see  \cite{SV}.
The function $\Phi$ is a function of variables 
$\bs  t = ( t_j^{(i)})$, where $ i = 1, \dots , r$, and $ j = 1, \dots , l_i$,
of variables $z = (z_1, \dots , z_n)$, weights $\bs \La$, and integers
$\bs l$. The main variables are $\bs t$, the other variables will be considered as 
parameters.

The  function $\Phi$ is symmetric with
respect to permutations of variables with the same upper index.

A point $\bs t$ with complex coordinates is called a {\it critical
  point} if the following system of algebraic equations is satisfied
\bean\label{Bethe eqn}
-\sum_{s=1}^n \frac{(\Lambda_s, \alpha_i)}{t_j^{(i)}-z_s}\ +\
\sum_{s,\ s\neq i}\sum_{k=1}^{l_s} \frac{(\alpha_s, \alpha_i)}{ t_j^{(i)} -t_k^{(s)}}\ +\
\sum_{s,\ s\neq j}\frac {(\alpha_i, \alpha_i)}{ t_j^{(i)} -t_s^{(i)}}
= 0, 
\eean
where $i = 1, \dots , r$ and $j = 1, \dots , l_i$.
In other words, the point $\bs t$ is a critical point if
\be
\left(\Phi^{-1}\frac{\partial \Phi }{\partial t_j^{(i)}}\right)(\bs t)=0, \qquad
i=1,\dots,r,\; {} \ {} j=1,\dots, l_i.
\ee

Note that the product of symmetric groups
$S_{\bs l}=S_{l_1}\times \dots \times S_{l_r}$ acts on the critical
set of the master function permuting the coordinates with the same
upper index. All orbits have the same cardinality $l_1! \cdots l_r!$\ .  
We do not make distinction between critical points in the same orbit.

\subsection{ Polynomials representing critical points }\label{PLCP}
Let $\bs t$ be a critical point of the master function 
$\Phi=\Phi(\bs{t}; \La_\infty)$.
Introduce an $r$-tuple of polynomials $\bs y=( y_1(x),$ $\dots ,$ $ y_r(x))$,
\bean\label{y}
y_i(x)\ =\ \prod_{j=1}^{l_i}(x-t_j^{(i)}).
\eean
Each polynomial is considered up to multiplication 
by a non-zero number.
The tuple defines a point in the direct product 
$\PCr$ of $r$ copies of the projective space associated with the vector 
space of polynomials in $x$. 
We say that the tuple $\bs y$ {\it represents the
critical point}.

It is convenient to think that the tuple $(1, \dots , 1)$ of constant polynomials
represents  in $\PCr$ the critical point of 
the master function with no variables. 
This corresponds to the case when $\bs l = (0, \dots , 0)$ and $\La_\infty =
\sum_{s=1}^n \La_s$.

Introduce polynomials
\bean\label{T}
{T}_i(x)=\prod_{s=1}^n(x-z_s)^{\langle \La_s, \al_i^\vee\rangle}, \qquad i = 1, \dots , r .
\eean

We say that a given tuple $\bs y$ is {\it generic with respect 
to weights  $\Lambda_1, \dots , \Lambda_n$  and points}
$z_1, \dots, z_n$ if
\begin{itemize} 
\item
each polynomial $y_i(x)$ has no multiple roots;
\item
all roots of $y_i(x)$ are different from roots of the polynomial $T_i$;
\item
any two polynomials $y_i(x)$, $y_j(x)$ have no common roots if $i\neq j$ and
$a_{i,j}\neq 0$.
\end{itemize}

A tuple  is generic if it represents  a critical point.

Let $W(f,g)=f'g-fg'$ be the Wronskian of functions $f, g$ of $x$.

A tuple $\bs y$  is called {\it fertile},
  if for every $i = 1, \dots , r$ there
exists a polynomial $\tilde y_i$ satisfying the equation
\bean\label{wronskian-critical eqn}
W( y_i, \tilde y_i)\ = \ T_i \ \prod_{j,\ j\neq i}y_j^{-\langle\al_j,\al_i^\vee\rangle}\ =\
\ T_i \ \prod_{j,\ j\neq i}y_j^{- a_{i,j}}\ .
\eean
The polynomial $\tilde y_i$ considered up to multiplication by a non-zero
number has the form
\bean\label{tilde}
\tilde y_i(x)\ =\ c_1\ y_i (x) \ \int \ T_i(x)\ \prod_{j=1}^r 
y_j(x)^{-a_{i,j}}
\ dx \ +\ c_2\ y_i(x)\ ,
\eean
where $c_1, c_2$ are complex  numbers, $c_1 \neq 0$.

If $\bs y$ is fertile and $i \in \{ 1, \dots, r\}$, then the tuple
\bean\label{simple}
\bs y^{[i]}\ =\ (y_1, \dots , \tilde y_i,\dots, y_r) \ 
\qquad \in \ \qquad \PCr 
\eean
is called {\it the immediate descendant} of $\bs y$
in the $i$-th direction.

\begin{thm}[\cite{MV1}]\label{fertile cor}
${}$

\begin{enumerate}
\item[(i)]
A generic tuple $\bs y = (y_1, \dots , y_r)$, with deg $y_j = l_j$, 
represents a critical point of the master function $\Phi (\bs t; \La_\infty)$,
with $\La_\infty = \sum_{s=1}^n \La_s - \sum_{j=1}^r \,l_j\, \al_j$,
if and only if it is fertile.
\item[(ii)] If $\bs y$ represents a critical point,  
$i \in \{ 1, \dots, r\}$, and $\bs y^{[i]}$ is 
an immediate descendant of $\bs y$, then $\bs y^{[i]}$ is fertile.
\end{enumerate}
\end{thm}

 Let $\bs y$ represent a critical point
of  $\Phi (\bs t; \La_\infty)$.
Let $i\in \{1, \dots , r\}$ and  let $\bs y^{[i]}$  be
an immediate descendant of $\bs y$ in the $i$-th direction. 
Denote $\tilde l_i = \deg \tilde y_i$
and $\La_\infty^{[i]} = \sum_{s=1}^n \La_s - \tilde l_i \al_i -
\sum_{j=1, \ {} j\neq i}^r l_j \al_j$. Assume that
$\bs y^{[i]}$ is generic, then 
$\bs y^{[i]}$ represents a critical point of the master  function
$\Phi (\bs t; \La_\infty^{[i]})$.
If $\tilde l_i\neq l_i$,  then 
$$
\La_\infty^{[i]}\  =\
s_i\cdot \La_\infty \ ,
$$
 where $s_i\cdot$ is the shifted action of the 
$i$-th reflection of the Weyl group.

\subsection{Simple reproduction procedure}
 Let $\bs y$ represent a critical point
of  $\Phi (\bs t, \La_\infty)$.
The tuples $\bs y^{[i]}\ =\ (y_1, \dots , \tilde y_i,\dots, y_r)\ {}\ \in \PCr$, 
where $\tilde y_i$ is given  by \Ref{tilde} and $c_1, c_2$ are arbitrary 
numbers, not both equal to zero,
form a one-parameter family. The parameter space of the family
is identified with the projective line $\bs P^1$ with projective coordinates 
$(c_1 : c_2)$.
We have a map
\bean\label{map}
Y_{\bs y, i}\ : \ \bs P^1 \ \to \PCr\ , 
\notag
\eean
which sends a point $c = (c_1 : c_2)$ to the corresponding tuple $\bs y^{[i]}$.
Almost all tuples $\bs y ^{[i]}$ are generic.
The exceptions form a finite set in $\bs P^1$.

Thus, starting with a tuple $\bs y$,
representing a critical point of the master function
$\Phi(\bs t; \La_\infty)$, and an index 
$i \in \{1, \dots , r\}$, we  construct  a  family 
$Y_{\bs y, i} : \bs P^1 \to \PCr$ of fertile tuples.
For almost all $c \in \bs P^1$ (with finitely many exceptions only), 
the tuple $Y_{\bs y, i}(c)$  represents a critical point of the master function
associated with integral dominant weights $\La_1, \dots , \La_n$, points
$z_1, \dots , z_n$, and a suitable weight at infinity.

We call this construction  {\it the
simple reproduction procedure in the $i$-th direction}.

\subsection{General reproduction procedure}\label{general procedure}

Assume that a tuple $\bs y \in \PCr$ represents a critical point of the master function
$\Phi(\bs t; \La_\infty)$. 

Let $\bs i = [i_1, i_2, \dots , i_k], \ i_j \in \{1, \dots , r\},$ 
be a sequence of natural numbers.
We define a $k$-parameter family of fertile tuples 
\bea\label{general map}
Y_{\bs y, \bs i} \ :\ (\bs P^1)^k \ \to \PCr
\eea
by induction on $k$,  starting at $\bs y$ and successively applying 
the simple reproduction procedure in directions
$i_1, \dots , i_k$. The image of this map is denoted 
$P_{\bs y, \bs i}$ .

For a given $\bs i = [i_1, \dots , i_k]$, almost all tuples 
$Y_{\bs y, \bs i} ( \bs c)$ represent  critical points of master functions associated to
weights $\La_1, \dots , \La_n$,  points $z_1, \dots , z_n$, and suitable weights at infinity.
Exceptional values of $c \in (\bs P^1)^k$ are contained in a proper algebraic subset.

It is easy to see that if $\bs i' = [i'_1, i'_2, \dots , i'_{k'}], \ 
 i_j \in \{1, \dots , r\},$  
is a sequence of natural numbers, 
and the sequence $\bs i'$ is contained in the sequence $\bs i$ as  an ordered subset, then
$P_{\bs y, \bs i'}$ is a subset of $P_{\bs y, \bs i}$.

The union 
\be 
P_{\bs y} \ = \ \cup_{\bs i} \ P_{\bs y, \bs i}\ \subset \PCr\ ,
\ee
where the summation is over all of sequences $\bs i$, 
is called {\it the population of critical points associated} with the Kac-Moody algebra $\g$, 
integral dominant
weights $\La_1, \dots , \La_n$, points $z_1, \dots , z_n$,  and {\it originated} at $y$.

If two populations intersect, then they coincide.

If the Weyl group is finite, then all tuples of a population consist of polynomials 
of bounded degree. Thus, if the Weyl group of $\g$ is finite, then a population
is a projective irreducible variety.

Every population $P$ has a tuple $\bs y = (y_1, \dots , y_r)$, $\deg y_i = l_i$, such
that the weight $\La_\infty  = \sum_{s=1}^n \La_s - \sum_{i=1}^r l_i \al_i$ is 
dominant integral, see \cite{MV1}.

\begin{conj}[\cite{MV1}]\label{CON}
Every population, associated with a Kac-Moody algebra $\g$, dominant integral weights
$\Lambda_1, \dots , \Lambda_n$,   points $z_1, \dots , z_n$,
 is an algebraic variety isomorphic to 
the flag variety associated to the Kac-Moody algebra
$^t\g$ which is Langlands dual to $\g$. Moreover, the parts of the family corresponding to 
tuples of polynomials with fixed degrees are isomorphic to 
Bruhat cells of the flag variety.
\end{conj}

The conjecture is proved  for the Lie algebras with root systems 
of types $A_r, B_r, C_r, G_2$ in \cite{MV1, BM}. In Theorems \ref{isomorphism}
and \ref{second main}
we prove this conjecture for any simple Lie algebra.

\subsection{Diagonal sequences of polynomials associated with a critical point 
and a sequence of indices}
We introduce notions which will be used
in Chapter \ref{solutions} to construct solutions of differential equations.

\begin{lemma}\label{generation}
Assume that a tuple $\bs y$ of non-zero polynomials
 represents a critical point of the master function
$\Phi(\bs t; \La_\infty)$. 
Let $\bs i = [i_1, i_2, \dots , i_k], \ 
 i_j \in \{1, \dots , r\},$  
be a sequence of natural numbers. Then there exist tuples
$\bs y^{[i_1]} = ( y_1^{[i_1]}, \dots , y_r^{[i_1]})$,
$\bs y^{[i_1, i_2]} = ( y_1^{[i_1,i_2]}, \dots , y_r^{[i_1,i_2]})$,
\dots ,
$\bs y^{[i_1, i_2, \dots , i_k]} = ( y_1^{[i_1,i_2, \dots , i_k]}, \dots , 
y_r^{[i_1,i_2, \dots , i_k]})$ in $\PCr$ such that
\begin{enumerate}
\item[(i)] 
\bea
% \label{ind eqn 1}
W( y_{i_1},  y_{i_1}^{[i_1]})\ = 
\ T_{i_1} \ \prod_{j,\ j\neq i_1} y_j^{- a_{i_1, j}}\ 
\eea
and $ y_j^{[i_1]} = y_j$ for $j\neq i_1$;
\item[(ii)] 
for $l = 2, \dots , k$, we have
\bea
% \label{ind eqn 2}
W( y_{i_l}^{[i_1, \dots , i_{l-1}]},  y_{i_l}^{[i_1, \dots , i_l]})\ = 
\ T_{i_l} \ \prod_{j,\ j\neq i_l}
(y_j^{[i_1, \dots , i_{l-1}]})^{-a_{i_l, j}}\ 
\eea
and $ y_{j}^{[i_1, \dots , i_{l}]} =  y_{j}^{[i_1, \dots , i_{l-1}]}$
for $j\neq i_l$.
\hfill $\square$
\end{enumerate}
\end{lemma}

The tuples 
$\bs y^{[i_1]}$,
$\bs y^{[i_1, i_2]}$,
\dots ,
$\bs y^{[i_1, i_2, \dots , i_k]}$
belong to the population  $P_{\bs y}$.
The tuple $\bs y^{[i_1]}$ is obtained from $\bs y$ by the $i_1$-th
simple generation procedure and  
for $l = 2, \dots , k$, the tuple $\bs y^{[i_1, \dots , i_{l}]}$
is obtained from $\bs y^{[i_1, \dots , i_{l-1}]}$
by the $i_l$-th simple generation procedure.

The sequence of tuples $\bs y^{[i_1]}$, $\bs y^{[i_1,i_2]}$, \dots ,
$\bs y^{[i_1, i_2, \dots , i_k]}$ satisfying Lemma \ref
{generation} will be called {\it associated with the critical point $\bs y$
and the sequence of indices $\bs i$}. The sequence of polynomials
$y_{i_1}^{[i_1]}$, $y_{i_2}^{[i_1,i_2]}$, \dots ,
$y_{i_k}^{[i_1, i_2, \dots , i_k]}$ 
 will be called {\it the diagonal sequence of polynomials associated with
the critical point $\bs y$ and the sequence of indices $\bs i$.} 
For a given $\bs y$ the diagonal sequence of
polynomials determine the sequence of tuples $\bs y^{[i_1]}$, $\bs y^{[i_1,i_2]}$, \dots ,
$\bs y^{[i_1, i_2, \dots , i_k]}$ uniquely.

There are many diagonal sequences of polynomials associated with a given critical point
and a given sequence of indices.

\section{Opers}\label{opers}

Let $\g=\g(A)$ be a Kac-Moody algebra with  simple roots
$\al_1, \dots , \al_r$ and simple coroots $\al_1^\vee, \dots , \al_r^\vee$.
Let $^t\g=\g({}^tA)$ be the Langlands dual algebra with
Chevalley generators $^t\h$, $E_1, \dots , E_r, F_1, \dots , F_r$, 
simple roots $^t\al_1, \dots , {}^t\al_r$ and simple coroots 
$^t\al_1^\vee, \dots , {}^t\al_r^\vee$.
Set $H_1\ =\ ^t\al_1^\vee, \ \dots ,\ H_r\ =\ {}^t\al_r^\vee$ and
\bea
I = F_1 + \dots + F_r , 
\qquad
\p = d/dx \ .
\eea
{\it A ${}^t\g$-oper} 
%associated with $^t\g$ 
is a differential operator of the form
\bea
D\ =\ \partial \ + \ I \ +\ V\ + \ W\ 
\eea
with $V \in M(^t\h)$ and $W\in \bar M({}^t\n_+)$. {\it A Miura ${}^t\g$-oper} 
%associated with $^t\g$ 
is a differential operator of the form
\bea
D\ =\ \partial \ + \ I \ +\ V\ 
\eea
with $V \in M(^t\h)$.

The differential operators of that type were
considered by V. Drinfeld and V. Sokolov in their study of the KdV type
equations \cite{DS}. On opers and Miura opers see \cite{BD, FFR, F1, F2, MV1, BM, MV4}.

For $u \in \bar M({}^t\n_+)$ and a ${}^t\g$-oper $D$, the differential operator
\bea
e^{\text{ad}\,u} \cdot D\ = \ D + [u,D] + \frac 12 [u,[u,D]] + \dots
\eea
is a ${}^t\g$-oper. The opers $D$ and $e^{\text{ad}\,u} \cdot D$ are called
{\it gauge equivalent}.
% of $D$ by $e^{\text{ad}\,u}$ is the oper $e^{\text{ad}\,u} D$.
%\bea
%e^{\text{ad}\,u} D\ = \ \p - \p(u) + e^{\text{ad}\,u} (I + V +   W)\ .
%\eea

Let $X$ be a $^t\g$-module with locally finite action of $^t\n_+$. 
Let $D$ be a ${}^t\g$-oper and $u\in \bar M({}^t\n_+)$. Then $D, \,
e^{\text{ad}\,u} \cdot D,\, e^{\pm\,u} $ determine linear operators on $M(X)$. Moreover, we have
\bea
e^{\text{ad}\,u} \cdot D\ = \ e^{u} \,D \,e^{-u} \ .
\eea

\begin{lemma} Let $D = \p + I + V$ be a Miura ${}^t\g$-oper. 
Let $g\in M(\C)$ and $i \in \{1, \dots , r\}$. Then
\bea
e^{\text{ad}\,( g E_i )} \cdot D \ =\ \partial \ +\ I \ + ( V\ +\ g\  H_i )\ -\
(g' \ + \ \langle {}^t\alpha_i , V \rangle g\ +\  g^2)\ E_i \ .
\eea
\end{lemma}
The proof is straightforward. 

\begin{corollary} 
Let $D = \p + I + V$ be a Miura ${}^t\g$-oper. Then the 
 ${}^t\g$-oper $e^{\text{ad}\,( g E_i )} \cdot D$ \linebreak
is a Miura oper 
if and only if the scalar rational function
$g$ satisfies the Ricatti equation
\bean\label{Ric}
g' + \langle {}^t\alpha_i , V \rangle g +  g^2 = 0 \ .
\eean
\end{corollary}

We say that the  Miura ${}^t\g$-oper $D$ is {\it deformable in the $i$-th direction}
if  equation \Ref{Ric} has a non-zero 
solution which is a rational function. 
%If $D$ is a
% Miura ${}^t\g$-oper deformable in the $i$-direction, and $g$ is a rational
%solution of \Ref{Ric}, then we say 
%that the  Miura ${}^t\g$-oper $e^{\text{ad}\,( g E_i )}  D$ is {\it the
%Miura deformation of $D$ in the $i$-th direction}.

\bigskip

Fix a collection of dominant integral weights  
$\bs\La=(\La_1, \dots , \La_n)$ of the Kac-Moody algebra $\g$,
and numbers $\bs z= \{z_1, \dots , z_n\} \subset \C$, \
$z_i\neq z_j$ if $i\neq j$.
Introduce polynomials 
$T_1(x), \dots , T_r(x)$ by formulas
\Ref{T}.

Let $\bs y = (y_1, \dots , y_r)$ be a tuple of non-zero polynomials.
We say that a Miura ${}^t\g$-oper $D = \partial  +  I  +  V$ is {\it associated with
weights $\bs\La$, numbers $\bs z$}, and the tuple
$\bs y = (y_1, \dots , y_r)$, 
if for every $i\in\{1, \dots , r\}$ we have
\bean\label{def of V}
\langle  \ {}^t\al_i, V \rangle\ =\ -\ \mathrm{log}'
\left( T_i \prod_{j=1}^r y_j^{-\langle\al_j,\al_i^\vee\rangle}
\right)\  =\ -\ \mathrm{log}'\left( T_i \prod_{j=1}^r y_j^{- a_{i,j}}
\right) ,
\eean 
cf. \Ref{wronskian-critical eqn}.

If a Miura oper $D$ is associated with weights $\bs\La$, 
numbers $\bs z$, and a generic tuple $\bs y = (y_1, \dots , y_r) \in 
\PCr$, then the tuple $\bs y$ is determined uniquely.
Indeed, the residues of the rational function 
$\langle  \ {}^t\al_i, V \rangle$ are positive exactly at the roots
of the polynomial $y_i$ and the residues are equal to the
multiplicities of roots of $y_i$ multiplied by two.

If $\g$ is a simple Lie algebra, then $D$ determines $\bs y$ uniquely even if
$\bs y$ is not generic. That fact follows from the invertibility of the Cartan matrix
of $\g$.

\begin{theorem}\label{main}
Let the  Miura ${}^t\g$-oper
$D = \partial  +  I  +  V$ be associated with
weights $\bs\La$, numbers $\bs z$, and the tuple $\bs y = (y_1, \dots , y_r)$.
Let $i\in \{1, \dots , r\}$. Then $D$ is deformable in the $i$-th direction
if and only if there exists a polynomial $\tilde y_i$ satisfying  
\Ref{wronskian-critical eqn}. Moreover, in that case any non-zero rational 
solution $g$ of the Ricatti equation \Ref{Ric} has the form
$g = \mathrm{log}' (\tilde y_i/ y_i)$ where $\tilde y_i$ is a solution of 
\Ref{wronskian-critical eqn}. If 
$g = \mathrm{log}' (\tilde y_i/ y_i)$, then
the Miura ${}^t\g$-oper
\bean\label{transformation}
e^{\text{ad}\,( g E_i )} \cdot D = \partial  + I  + ( V \,+\, g\  H_i )
\eean
is associated with weights $\bs\La$, numbers $\bs z$,
and the tuple $\bs y^{[i]} = (y_1, \dots ,\tilde y_i,$
$ \dots , y_r)$, where the tuple $\bs y^{[i]}$ is called in Section
\ref{PLCP} an immediate descendant of $\bs y$ in the $i$-th direction, see \Ref{simple}.

\end{theorem}

\begin{proof} Write \Ref{Ric} as
\bean\label{Ric1}
g'/g  +  g  =  \mathrm{log}'\left( T_i \prod_{j=1}^r y_j^{- a_{i,j}}
\right) \ .
\eean
If $g$ is a rational function, then $g \to 0$ as $x \to \infty$ and all poles 
of $g$ are simple. Moreover, the residue of $g$ at any point is an integer. Hence
$g = c'/c$ for a suitable rational function $c$. Then 
\bean\label{Ric2}
c\ = \ \int T_i(x) \prod_{j=1}^r y_j(x)^{-a_{i,j}} dx \ 
\eean
and equation \Ref{wronskian-critical eqn} has a polynomial solution 
$\tilde y_i = - c y_i$. Conversely if  equation
\Ref{wronskian-critical eqn} has a polynomial solution $\tilde y_i$, then
the function $c$ in \Ref{Ric2} is rational. Then $g = c'/c$ is a rational solution of 
equation \Ref{Ric}.

Let $g = \mathrm{log}' (c) = \mathrm{log}' (\tilde y_i/ y_i)$, 
where $\tilde y_i$ is a solution of \Ref{wronskian-critical eqn}. Then
\bea
e^{\text{ad}\,( g E_i )} \cdot D = \partial  + I  + 
( V \,+\, \mathrm{log}'(\tilde y_i / y_i)\ H_i )
\eea
and
\bea
\langle \  {}^t\al_k, V \rangle \ + \ \mathrm{log}'(\tilde y_i / y_i)\,
\langle \  {}^t\al_k,\ {}^t\alpha_i^\vee  \rangle \
 = \ -\ \mathrm{log}'\left(
 T_k \,\tilde y_i^{-\langle\al_i,\al_k^\vee\rangle}\,
\prod_{j=1, \ j\neq i }^r y_j^{-\langle\al_j,\al_k^\vee\rangle}
\right) .
\eea
\end{proof}

Note that if equation \Ref{Ric} has one non-zero rational solution
$g = c'/c$ with rational $c$, 
then other non-zero (rational) solutions have the form
$g = c'/(c + \text{const})$.

\begin{corollary}
Let the  Miura ${}^t\g$-oper
$D = \partial  +  I  +  V$ be associated with
weights $\bs\La$, numbers $\bs z$, and the tuple $\bs y = (y_1, \dots , y_r)$.
Then $D$ is deformable in all directions from 1 to $r$
if and only if the tuple $\bs y$ is fertile.
\end{corollary}

\begin{corollary}
Let the  Miura ${}^t\g$-oper
$D = \partial  +  I  +  V$ be associated with
weights $\bs\La$, numbers $\bs z$, and the tuple $\bs y = (y_1, \dots , y_r)$.
Let the  tuple $\bs y = (y_1, \dots , y_r)$ be generic in the sense of
Section \ref{PLCP}. Then
 $D$ is deformable in all directions from 1 to $r$
if and only if the tuple $\bs y$
represents  a critical point of the master function
\Ref{master} associated with parameters $\bs z, \bs \La, \La_\infty$.
\end{corollary}

Let the  Miura ${}^t\g$-oper
$D = \partial  +  I  +  V$ be associated with
weights $\bs\La$, numbers $\bs z$, and the tuple $\bs y = (y_1, \dots , y_r)$.
Let the  tuple $\bs y = (y_1, \dots , y_r)$ 
represent  a critical point of the master function
\Ref{master} associated with parameters $\bs z, \bs \La, \La_\infty$.
Let $\om^0 $ be the variety 
of all Miura opers which can be obtained from $D$ by a sequence of deformations
in directions $i_1, \dots , i_N$ where $N$ is any positive integer and all
$i_j$ lie in $\{1, \dots , r\}$.

\begin{corollary}\label{om0}
For a simple Lie algebra $\g$ 
the variety $\om^0$ is isomorphic to the population of critical points 
originated at $\bs y$.
\end{corollary}

\section{Miura opers and flag varieties}

In this section we assume that $\g$ and ${}^t\g$ are simple Lie algebras although
most of considerations can be extended to  Kac-Moody algebras.

Let $\tG$ be the complex simply connected
 Lie group with Lie algebra $\tg$. Let ${}^tB_\pm,$ $ {}^tN_\pm,$ $ 
{}^tH$ be 
the subgroups with Lie algebras ${}^t\B_\pm, {}^t\n_\pm, {}^t\h$, respectively.

\subsection{Triviality of the monodromy}\label{monodromy}
Let $D = \p + I + V$ be a Miura $\tg$-oper.
% associated with
%weights $\bs\La$, numbers $\bs z$, and a tuple $\bs y = (y_1, \dots , y_r)$. 
%Assume that the Miura oper is deformable in all directions from 1 to $r$.
Let $\bP$ be the complex projective line. Consider $D$ as a
$\tG$-connection $\nabla_D$ on the trivial principal $\tG$-bundle $p :
\tG \times \bP \to  \bP$.  The connection has singularities at
the set $\mathrm{Sing} \subset \bP$ where the function $V$ has poles. 
Choose a regular point $x_0 \in \bP - \infty$ of the connection. 
Parallel translations with respect
to the connection
define the monodromy representation $\pi (\bP - \mathrm{Sing}) \to {}^tG$.
Its image is called {\it the monodromy group}.

\begin{theorem}\label{Monodromy}
Assume that the Miura $\tg$-oper $D$ is  associated with
weights $\bs\La$, numbers $\bs z$, and a tuple $\bs y = (y_1, \dots , y_r)$
as in Section \ref{opers}. 
Assume that the tuple $\bs y = (y_1, \dots , y_r)$ is generic in the sense of
Section \ref{PLCP}. 
Assume that the Miura oper $D$ is deformable in all directions from 1 to $r$.
Then the monodromy group of $\nabla_D$ belongs to the center of ${}^tG$.

\end{theorem}

\begin{proof} It is known that the intersection of all 
of the Borel subgroups in $\tG$ is
the center of $\tG$, see \cite{B, H}. 
We show that the monodromy of $\nabla_D$ lies in the intersection
of all of the Borel subgroups.

Let $\id \in \tG$ be the identity element.
Let $\bar Y(x) \in \tG$ be the (possibly multi-valued) solution of the equation
$D Y = (\p + I + V)Y = 0$ such that $Y(x_0)= \id$. 
Since $D$ is a Miura oper we have for any regular $x$ the equality of sets
 $\bar Y(x)\,\tB = \tB$.
Hence if $m\in \tG$ is an element of the monodromy group of $\nabla_D$, then
$m\ \tB = \tB$ and hence $m\in \tB$.

Let $i\in\{1, \dots , r\}$ and let $g_i\in M(\C)$ be a solution of the Ricatti equation
\Ref{Ric}. Assume that $g_i$ is regular at $x_0$. Then the $\tG$-valued function
\bea
%Y^{[i]}_{g_i}(x)\ =\
 e^{ g_i(x) E_i }\ \bar Y(x)\ e^{ - g_i(x_0) E_i } 
\eea
is the 
solution of the equation $(e^{\text{ad}\,( g_i(x) E_i )}\cdot D)\, Y  =  0$ 
such that $Y(x_0)= \id$. Since $e^{\text{ad}\,( g_i(x) E_i )} \cdot D$
is a Miura oper we have an equality of sets
$Y^{[i]}_{g_i}(x)\,\tB = \tB$ for any  $x$ at which 
$e^{\text{ad}\,( g_i(x) E_i )}\cdot D$ is regular.
Hence 
for any element $m$ of the monodromy group of $\nabla_D$
we have
\bea
e^{ g_i(x_0) E_i }\ m \ e^{ - g_i(x_0) E_i }\ ^tB_- \ {} \subset \ ^tB_-
\eea
or $m \in\ {} e^{ - g_i(x_0) E_i }\ ^tB_-\ e^{ g_i(x_0) E_i }$.

Now consider the Miura oper
\bea
e^{\text{ad}\,( g_i(x) E_i )} \cdot D = \p + I + V_{i;  g_i}
\eea
where $V_{i;  g_i}$ is the ${}^t\h$-part of $e^{\text{ad}\,( g_i(x) E_i )}\cdot D$.
Let $j\in\{1, \dots , r\}$ and let $g_{i,j; g_i}\in M(\C)$ be a solution of the $j$-th
Ricatti equation 
\bea
g' + \langle {}^t\al_i, V_{i;  g_i}\rangle g + g^2 = 0
\eea
associated with the Miura oper $e^{\text{ad}\,( g_i(x) E_i )} \cdot D$. 
Assume that $g_{i,j; g_i}$ is regular at $x_0$. Then the $\tG$-valued function
\bea
%Y^{[j]}_{i,j;\ g_i, g_{i,j; g_i}}(x)\ =\
 e^{ g_{i,j; g_i} (x) E_j } e^{ g_i(x) E_i }\ \bar Y(x)\ e^{ - g_i(x_0) E_i } 
 e^{ -g_{i,j; g_i} (x_0) E_j }
\eea
 is the 
solution of the equation 
$(e^{\text{ad}\,(  g_{i,j; g_i}(x) E_j )}e^{\text{ad}\,( g_i(x) E_i )}\cdot D)\, Y  =  0$ 
such that $Y(x_0)= \id$.
Repeating the previous argument we conclude that any element $m$
of the monodromy group of $\nabla_D$ lies in the Borel subgroup
\bea
e^{ - g_i(x_0) E_i }  e^{ -g_{i,j; g_i} (x_0) E_j } \ \tB
\ e^{ g_{i,j; g_i} (x_0) E_j } e^{ g_i(x_0) E_i } .
\eea

Every $u\in \ {}^tN_+$ is a product of elements of the form 
$e^{c_i E_i }$ for $i \in \{1, \dots , r\}$ and $c_i\in\C$.
Every $c_i$ can be taken as the initial condition for a solution of the suitable
$i$-th Ricatti equation. Therefore 
the iteration of the previous reason shows that every element 
of the monodromy group of $\nabla_D$ lies in every Borel subgroup
of the form $u^{-1} (\, \tB) u$, $u\in \ {}^tN_+$.
The Borel subgroups in $\tG$ of the form $u^{-1}(\,\tB ) u$,
\ $u\in \ {}^tN_+$, form an open dense subset in the flag variety
of all of the Borel subgroups. Hence the monodromy lies in the intersection of all
of the Borel subgroups.

\end{proof}

\subsection{Gauge equivalent Miura opers}\label{gauge}
As in Section \ref{monodromy}, let $D$ be the 
Miura $\tG$-oper   associated with
weights $\bs\La$, numbers $\bs z$, and a tuple $\bs y = (y_1, \dots , y_r)$.
We assume that the tuple $\bs y = (y_1, \dots , y_r)$ is generic in the sense of
Section \ref{PLCP} and 
the Miura oper $D$ is deformable in all directions from 1 to $r$.

Consider the variety $\om$ of all Miura opers gauge equivalent to $D$. 
If $D' \in \om$, then there exists a rational $N_+$-valued function
$v$ on $\bP$ such that $D' \ =\ v\, D\, v^{-1}$. In that case we denote $D'$ by
$D^v$. 

Let $\om^0 \subseteq \om$ be the subvariety 
of all Miura opers which can be obtained from $D$ by a sequence of deformations
in directions $i_1, \dots , i_N$ where $N$ is a non-negative integer and all
$i_j$ lie in $\{1, \dots , r\}$. By Corollary \ref{om0} the subvariety
$\om^0$ is isomorphic to the population of critical points originated at $\bs y$.

\medskip

The connection $\nabla_D$ is regular at $x_0\in \bP$ if $x_0$ does not lie in
$\{z_1, \dots , z_n, \infty\}$ and $x_0$ is not a root of some of polynomials
$y_1, \dots , y_r$.

Consider the trivial bundle  $p' : (\tG/\tB) \times \mathbb P^1 \to \bs \bP$ associated
with the bundle $p$. The fiber of $p'$ is
the flag variety $\tG/\tB$. The connection $\nabla_D$ induces
a connection $\nabla_D'$ on $p'$. The monodromy of $\nabla_D'$ is trivial by Theorem
\ref{Monodromy}.
Thus the variety $\Gamma$ of global horizontal sections of $\nabla_D'$ is identified with
the fiber $(p')^{-1}(x_0)$ over any $x_0$ which is a regular point of 
the connection. Thus $\Gamma$ is isomorphic to $\tG/\tB$.

Any $\tG$-valued rational function $v$ defines a section
\bean\label{section}
S_v \ :\ x\ \mapsto\ v(x)^{-1}\ \tB \ \times\ x
\eean
 of $p'$ over the set of regular points of $v$.
The section $S_v$ is also well defined over the poles of $v$ since $\tG/\tB$ 
is a projective variety.

If $D^v \in \om$, then the section $S_v$ is horizontal with respect to
$\nabla_D'$, cf. the proof of Theorem \ref{Monodromy}.
Thus we have a map 
\bea
S : \om \to \Gamma , 
\qquad
D^v \mapsto S_v \ .
\eea

\begin{theorem}\label{isomorphism}
The map $ S : \om \to \Gamma $ is an isomorphism and $\om^0=\om$.
\end{theorem}

\begin{proof}
Let $D^{v_1}, D^{v_2} \in \om$. Assume that the images of
$D^{v_1}$ and $D^{v_2}$ under the map $S$ coincide.
Assume that $v_1, v_2, D$ are regular
at $x_0\in \bP$. The equality $S_{v_1}(x_0)= S_{v_2}(x_0)$
means that  $v_1(x_0)^{-1}\,\tB = v_2(x_0)^{-1}\,\tB$. Then
$v_1(x_0) = v_2(x_0)$. Hence $v_1=v_2$ and $D^{v_1} = D^{v_2}$.
That proves the injectivity of $S$.

Let $x_0$ be a regular point of $D$ in $\bP - \infty$. For any $u \in \ \tN$
there exists a rational $\tN$-valued function $v$ such that $v(x_0) = u$
and $D^v \in \om$ and $D^v$ is obtained from $D$
by a sequence of deformations in some directions $i_1, \dots , i_N$,
see the proof of Theorem \ref{Monodromy}. Thus the set
\bea
Im (x_0) \  =\ \{ S_v(x_0) \in \ (\tG/\tB)\times x_0 \ | \ D^v \in \om^0 \}
\eea
contains the set $((\tN \ \tB)/\tB )\times x_0\ \subset\ (\tG/\tB)\times x_0$.
The set $Im(x_0)$ is  the image with respect to $S$
of a population of critical points. Hence it is closed
as the image of a closed variety. On the other hand
the set $((\tN \ \tB )/\tB)\times x_0$ is dense in $(\tG/\tB)\times x_0$. Hence
$Im(x_0) = (\tG/\tB)\times x_0$ and
$S(\om^0) = \Gamma$. Therefore 
 $\om^0 = \om$ since the map $S$ is injective.
\end{proof}

\subsection{Remarks on the isomorphism}\label{REMARKS}
Let $\g$ be a simple Lie algebra. Let $P_{\bs y_0}$ be the population of critical points
originated at a tuple $\bs y^0$. We assume that the tuple $\bs y^0$ is generic in the sense of
Section \ref{PLCP}. Theorem \ref{isomorphism} says that the population 
$P_{\bs y^0}$ is isomorphic to the flag variety
$\tG/\tB$. The isomorphism is constructed in three steps. 
If $\bs y' \in P_{\bs y}$ is a point of the population, then one assigns to it
the associated Miura oper $D_{\bs y'}$ as in Section \ref{opers}, see also Lemma
\ref{reduced}. By Theorem \ref{main} we have
$D_{\bs y'}\ =\ v\ D_{\bs y^0}\ v^{-1}$ for a suitable rational function $v : \bP \to \ {}^tN_+$.
To the Miura oper $D^v$ one assigns the section $S_v \in \Gamma$ by formula \Ref{section}.
Then one chooses a point $x_0\in \C$,
regular with respect to the connection $\nabla_{D_{\bs y^0}}$,
and assigns to a section $S^v$ its value $S^v(x_0) \in (\tG/\tB) \times x_0$.
The resulting composition 
\bea
\phi_{\bs y^0, x_0}\ :\ P_{\bs y} \ \to \ \tG/\tB 
\eea
is an isomorphism according to Theorem \ref{isomorphism}.

\begin{lemma}\label{dep 1}
If $x_0, x_1 \in \C$ are points regular with respect to 
$\nabla_{D_{\bs y^0}}$,
then there exists an element $g\in \tB$ such that \ {}\
$\phi_{\bs y^0, x_1}\ = \ g \ \phi_{\bs y^0,  x_0}$.
\end{lemma}

\begin{proof} Let $Y$ be the $\tG$-valued solution of the equation
$D_{\bs y^0} Y = 0$ such that $Y(x_0) = \id$. Then $Y(x) \in \tB$ 
for all  $x$. 
If $D^v \in \om_{\bs y^0}$, then $S^v$ is a horizontal section
of $\nabla_{D_{\bs y^0}}'$. Thus it has the form
$ x \ \mapsto\ (Y(x)\,u \ \tB)\times x$ for a suitable element $u\in \tG$.
Hence $\phi_{\bs y^0, x_0} ( \bs y' ) \ =\  Y(x_0)\,u\ \tB$
and
$\phi_{\bs y^0, x_1} ( \bs y' )\ =\ Y(x_1)\,u\ \tB$.
We conclude that
$\phi_{\bs y^0, x_1}\ =\ Y(x_1)Y(x_0)^{-1} \,\phi_{\bs y^0, x_0}$.
\end{proof}

Let $\bs y^1 $ be a point of $P_{\bs y^0}$. 
Let $P_{\bs y^1}$ be the population originated at $\bs y^1$. 
 We have  $P_{\bs y^0} =  P_{\bs y^1}$. 

\begin{lemma}\label{dep 2}
Let $x_0 \in \C$ be regular with respect to both connections 
$\nabla_{D_{\bs y^0}}$ and $\nabla_{D_{\bs y^1}}$. Then there
exists an element $g \in \ {}^tB_+$ such that
$\phi_{\bs y^1, x_0}\ =\ g\ \phi_{\bs y^0, x_0}$.
\end{lemma}

\begin{proof}
We have $D_{\bs y^1}\ =\ w \ D_{\bs y^0}\ w^{-1}$
for a suitable rational function $w\ :\ \bP\ \to \ {}^tN_+$.
If $Y_0(x)$ is the $\tG$-valued solution of the equation
$D_{\bs y^0} Y = 0$ such that $Y_0(x_0) = \id$, then 
$Y_1(x)\ =\ w(x)\,Y(x)\, w(x_0)^{-1}$ is the 
$\tG$-valued solution of the equation
$D_{\bs y^1} Y = 0$ such that $Y_1(x_0) = \id$. 

Let $\bs y' \in P_{\bs y^0}$ and
$D_{\bs y'}\ =\ v\ D_{\bs y^0}\ v^{-1}$
for a suitable rational function $v\ :\ \bP\ \to \ {}^tN_+$.
Then $D_{\bs y'}\ =\ vw^{-1}\ D_{\bs y^1}\ wv^{-1}$.
Hence
$\phi_{\bs y^0, x_0} ( \bs y' )\ =\ v(x_0)^{-1}\,\tB$
and
$\phi_{\bs y^1, x_0} ( \bs y' )\ =\ w(x_0)v(x_0)^{-1}\,\tB$.
Therefore, 
$\phi_{\bs y^1, x_0} \ =\ w(x_0)\, \phi_{\bs y^0, x_0}$.
\end{proof}

\section{Bruhat Cells}\label{cells}

\subsection{Properties of Bruhat cells} Let $\g$ be a simple Lie algebra. 
For an element $w$ of the Weyl group $W$, the set
\bea
B_w\ = \ \tB \, w\, \tB \ {} \subset \ {} \tG/\tB
\eea
is called {\it the Bruhat cell} associated to $w$. The Bruhat cells
form a cell decomposition of the flag variety $\tG/\tB$.

For $w\in W$ denote $l(w)$ the length of $w$. We have 
$\dim B_w = l(w)$.

Let $s_1, \dots , s_r \in W$ be the generating reflections of the Weyl group.

For $v \in \tG/\tB$ and $i \in \{1, \dots , r\}$ consider the rational
curve 
\bea 
\C\ \to\ \tG/\tB, 
\qquad
 c \ \mapsto \  e^{ c E_i}\,v \ .
\eea
 The limit of
$e^{ c E_i}\,v$ is well defined as $c \to \infty$, since
$\tG/\tB$ is a projective variety.

We  need the following standard property of Bruhat cells.

\begin{lemma}\label{BRUHAT}
Let $s_i, w \in W$ be such that $l( s_i w ) = l(w) + 1$. Then
\bea
B_{ s_i w }\ =\
\{\ e^{ c E_i }\,v\ | \ v \in B_w,\ c\in \{\bP - 0\}\ \} \ .
\eea
\hfill
$\square$
\end{lemma}

\begin{corollary}\label{Bruhat}
Let $w = s_{i_1} \cdots s_{i_k}$ be a reduced decomposition of $w\in W$. Then
\bea
B_{ w } =
\{ \lim_{c_1 \to c^0_1}\,\dots\,\lim_{c_k \to c^0_k}
e^{ c_1 E_{i_1} }\,\cdots \,e^{ c_k E_{i_k} }\,\tB \ \in \ \tG/\tB \  | \  c^0_1, \dots , 
c^0_k
\in \{\bP - 0\} \}.
\eea
\end{corollary}

Introduce the map
\bea
f_{i_1, \dots , i_k} :  (\C- 0)^k \to B_{s_{i_1}\cdots\, s_{i_k}} ,
\qquad
(c_1, \dots , c_k) \mapsto
e^{ c_1 E_{i_1} }\,\cdots \,e^{ c_k E_{i_k} }\  \tB \ .
\eea

\subsection{Populations and Bruhat cells}
Let $P$ be a population of critical points
associated with weights $\bs\La$, numbers $\bs z$.
Let $T_1, \dots , T_r$ be the polynomials defined by \Ref{T}.
Let $\bs y^0 = (y^0_1, \dots , y^0_r) \in P$ with
 $l_i = \deg\, y^0_i$ for  $i \in \{1, \dots , r\}$.
Assume that the weight at infinity of $\bs y^0$,
\bea
\La_\infty \ =\ \sum_{i=1}^n \La_i\ - \
\sum_{i=1}^r\ l_i\,\al_i\ ,
\eea
is integral dominant, see Section \ref{crit pts}. Such $\bs y^0$ exists according to
\cite{MV1}.
For $w\in W$ consider the weight $w \cdot \La_\infty$, where
$w \cdot $ is the shifted action of $w$ on $\h^*$. Write
\bea
w \cdot \La_\infty \ =\ \sum_{i=1}^n \La_i\ - \
\sum_{i=1}^r\ l^w_i\,\al_i .
\eea
Set
\bea
P_w\ =\ \{\ \bs y = (y_1, \dots , y_r)\,\in \,P\ | \ 
\deg y_i = l^w_i, \ i = 1, \dots , r \ \}\ .
\eea

Consider the trivial bundle  $p' : (\tG/\tB) \times \mathbb P^1 \to \bs \bP$ 
with  connection $\nabla_{D_{\bs y^0}}'$. Consider the Bruhat cell decomposition
of fibers of $p'$. 

Let $x_0\in \C$ be such that $T_i(x_0)\neq 0$ and $y^0_i(x_0)\neq 0$ for 
$i = 1, \dots , r$. The point
$x_0 \in \C$ is a regular point of the connection $\nabla_{D_{\bs y^0}}'$.
Let 
\bea
\phi_{\bs y^0, x_0}\ :\ P\ \to \ \tG/\tB
\eea
be the isomorphism defined in Section \ref{REMARKS}.

\begin{theorem}\label{second main}
For $w \in W$ we have 
\bea
\phi_{\bs y^0,\, x_0} ( P_w )\ =\ B_{w^{-1}}\ .
\eea
\end{theorem}

\begin{corollary}\label{main corollary}
Let $\La_1, \dots , \La_n, \La_\infty$ be integral dominant $\g$-weights.
Let $z_1, \dots , z_n$ be distinct complex numbers. Let $w\in W$. 
Consider the master function $\Phi(\bs t; \bs z, \bs \La, w\cdot \La_\infty)$.
Let $K$ be a connected component of the critical set of the master function.
For each $\bs t \in K$ consider the tuple $\bs y_{\bs t} \in (\C[x])^r$ 
of monic polynomials
representing the critical point $\bs t$. Then the closure of the
set $\{\ \bs y_{\bs t} \ | \ \bs t \in K \ \}$
is an $l(w)$-dimensional cell.
\end{corollary}

\subsection{Proof of Theorem \ref{second main}}

\begin{lemma}\label{invariance}
For $w \in W$,  the subset $B_w \times \bP \ \subset \
(\tG/\tB) \times \bP$ is invariant with respect to
the connection $\nabla_{D_{\bs y^0}}'$.
\end{lemma}

\begin{proof} Let $Y$ be the $\tG$-valued solution of the equation
$D_{\bs y^0} Y = 0$ such that $Y(x_0) = \id$.  Then $Y(x) \in \tB$ 
for all  $x$. The horizontal sections of 
$\nabla_{D_{\bs y^0}}'$ have the form
$ x \ \mapsto\ (Y(x)\,u \ \tB)\times x$ for a suitable element $u\in \tG$.
Hence if $u \ \tB \ \in B_w$, then
$Y(x)u \ \tB\  \in B_w$ for all $x$.
\end{proof}

Let $w = s_{i_k} \cdots s_{i_1}$ be a reduced decomposition of $w\in W$. 
For $b = 1, \dots , k$ set
\bea
(s_{i_b} \cdots s_{i_1}) \cdot \La_\infty \ =\
\sum_{i=1}^n \La_i\ - \ \sum_{i=1}^r\ l^b_i\,\al_i .
\eea
From \cite{BGG} it follows that $l^{1}_{i_1} > l_{i_1}$ and
$l^{b}_{i_b} > l^{b-1}_{i_b}$ for $b = 2, \dots , k$.

For $b = 1, \dots , k$ define by induction on $b$
a family of tuples of polynomials  depending on
complex parameters $c_1, \dots , c_b$. 
Namely, let $\tilde y_{i_1}$ be a polynomial satisfying equation
\bea
W(\, y^0_{i_1}, \tilde y_{i_1}\,)\ = 
\ T_{i_1}\ \prod_{j,\ j\neq i_1} ( \,y^0_j\, )^{- a_{i_1,j}}\ .
\eea
We fix $\tilde y_{i_1}$ assuming that the coefficient of $x^{l_{i_1}}$ in $\tilde y_{i_1}$
is equal to zero.
Set
$\bs y^{1;\, c_1}  = (y^{1;\, c_1}_1, \dots , y^{1;\, c_1}_r)$, where
\bea
y^{1;\, c_1}_{i_1}(x)\ =\ \tilde y_{i_1}(x) \ +\ c_1\, y_{i_1}^0(x)
\qquad
\text{and}\
y^{1;\, c_1}_{j}(x)\ =\ y_{j}^0(x)
\ {} \text{for}\ {}
j \neq i_1 \ .
\eea
%The family $\bs y^{1;\, c_1}$ does not depend on the choice of $\tilde y_{i_1}$.

Assume that the family $\bs y^{b-1;\, c_1, \dots , c_{b-1}}$ is already defined.
Let $\tilde y_{i_b}^{\ b-1;\, c_1, \dots , c_{b-1}}$ 
be a polynomial satisfying equation
\bea
W(\, y_{i_b}^{b-1;\, c_1, \dots , c_{b-1}}, \tilde y_{i_b}^{\ b-1;\, c_1, \dots , c_{b-1}}\,)\
 = 
\ T_{i_b} \ \prod_{j,\ j\neq i_b} (\,y_{j}^{b-1;\, c_1, \dots , c_{b-1}}\,)^{- a_{i_b,j}}\ .
\eea
We fix $\tilde y_{i_d}^{\ d-1;\, c_1, \dots , c_{d-1}}$
 assuming that the coefficient of $x^{l^{d-1}_{i_d}}$ in 
$\tilde y_{i_d}^{\ d-1;\, c_1, \dots , c_{d-1}}$
is equal to zero.
Set $\bs y^{b;\, c_1, \dots , c_{b}}  = (y^{b;\, c_1, \dots , c_b}_1, 
\dots , y^{b;\, c_1, \dots , c_b}_r)$, where
\bea
y^{b;\, c_1, \dots , c_{b}}_{i_b}(x)\ =\ \tilde y_{i_b}^{\ b-1;\, c_1, \dots , c_{b}}
(x) \ +\ c_b\,
y_{i_b}^{b-1;\, c_1, \dots , c_{b-1}}(x)
\eea
and
\bea
y^{b; \,c_1, \dots , c_{b}}_{j}(x)\ =\  
y_{j}^{b-1;\, c_1, \dots , c_{b-1}}(x)
\ {}\ \text{for} \ {}\ 
j \neq i_b \ .
\eea
%The family  $\bs y^{b;\, c_1, \dots , c_{b}}$
%does not depend on the choice of 
%$\tilde y_{i_b}^{\ b-1;\, c_1, \dots , c_{b}}$.

The $b$-th family  is obtained from the $(b-1)$-st family by the generation procedure in
the $i_b$-th direction, see Section \ref{general procedure}. For any $c_1, \dots , c_k$ the tuple
$\bs y^{k;\, c_1, \dots , c_{k}}$ lies in $P$.

For any $c_1, \dots , c_k$ and any $i \in \{1, \dots , r\}$, we have
\bea
\deg\
y^{k; \,c_1, \dots , c_{k}}_{i}(x)\ =\  l^w_i\ .
\eea

Set
\bea
P^{[i_1, \dots , i_k]}\ = \
\{ \ \bs y^{k;\, c_1, \dots , c_{k}} \ |\
c_1, \dots , c_k \in \C\ \} \ .
\eea

\begin{proposition}\label{MAIN}
 We have
\bea
\phi_{\bs y^0, x_0} ( P^{[i_1, \dots , i_k]} )\ =\ B_{w^{-1}}\ .
\eea
\end{proposition}

\noindent
{\it Proof of Proposition \ref{MAIN}.}
Let $D_{\bs y^{k;\, c_1, \dots , c_{k}}}$ be the Miura oper associated with
the tuple $\bs y^{k;\, c_1, \dots , c_{k}}$, then
\bean\label{FORMULA}
{}
\\
D_{\bs y^{k;\, c_1, \dots , c_{k}}} =
\exp\!\left( \mathrm{ad}\,  \mathrm{log}'  \left(\frac{ y_{i_k}^{k;\, c_1, \dots , c_{k}}}
{ y_{i_k}^{k-1;\, c_1, \dots , c_{k-1}}}\right) E_{i_k} \right)
\dots\,
\phantom{aaaaaaaaaaaaaaaaaaaaaaaa}%
\notag
\eean
\bea
\phantom{aaaa}
\exp\!\left( \mathrm{ad}\,
    \mathrm{log}'  \left( \frac{ y_{i_2}^{2;\, c_1,c_2} }{ y_{i_2}^{1;\, c_1}}
\right) E_{i_2} \right) 
\exp\!\left(  \mathrm{ad}\,  \mathrm{log}'  \left(
 \frac{ y_{i_1}^{1;\,c_1}}{ y^0_{i_1}}\right) E_{i_1} \right)\ \cdot\ 
D_{\bs y^0} \ =
\eea
\bea
\exp\!\left(- \mathrm{ad}   \left(\frac{ 
\ T_{i_k} \ \prod_{j,\ j\neq i_k} (y_{j}^{k-1;\, c_1, \dots , c_{k-1}})^{- a_{i_k,j}}}
{y_{i_k}^{k;\, c_1, \dots , c_{k}}\, y_{i_k}^{k-1;\, c_1, \dots , c_{k-1}}}\right) E_{i_k} \right)
\dots\,
\phantom{aaaaaaaaaaaaaaaaaaaaaaaaaaaaaaaaaaaaaaaaaaaaaaaaaaaaaaaaa}
\eea
\bea
%\phantom{aaaaaaaaaaaaaaaaaaaaaa}
\exp\!\left(\!-\mathrm{ad}\!\left(\frac{ 
T_{i_2} \ \prod_{j,\ j\neq i_2} (y_{j}^{1;\, c_1})^{- a_{i_2,j}}}
{y_{i_2}^{2;\, c_1,c_2}  y_{i_2}^{1;\, c_1}}
\right)\!E_{i_2}\!\!\right) 
%\cdot
%\phantom{aaaaaaaaaaaaaaaaaaaaaa}
%\eea
%\bea
%\phantom{aaaaaaaaaaaaaaaaaaaaaa}
\exp\!\left(\!-\mathrm{ad}\!\left(\frac{ T_{i_1} \prod_{j,\ j\neq i_1} (y^0_j)^{- a_{i_1,j}}}
{y_{i_1}^{1;\,c_1} y^0_{i_1}}\right)\!E_{i_1}\!\!\right)\!\cdot\! 
D_{\bs y^0},
\eea
see Theorem \ref{main}.

Introduce the rational map 
\bea
g : \C^{k+1} \to \C^k, \ {}\ {} \ (x; c_1, \dots , c_k) \mapsto
(g_1(x; c_1), \dots , g_k(x; c_1, \dots , c_k))
\eea
 where
\bea
g_1(x; c_1)
 =  \frac{ T_{i_1}( x ) \prod_{j,\ j\neq i_1} (\,y^0_j(x)\,)^{- a_{i_1,j}}}
{y_{i_1}^{1;\,c_1}(x)\ y^0_{i_1}(x)} ,
\eea
\bea
g_b(x; c_1, \dots , c_b) = 
\frac{ 
\ T_{i_b}(x) \ \prod_{j,\ j\neq i_b} (\,y_{j}^{b-1;\, c_1, \dots , c_{b-1}}(x)\,)^{- a_{i_b,j}}}
{y_{i_b}^{b;\, c_1, \dots , c_{b}}(x)\  y_{i_b}^{b-1;\, c_1, \dots , c_{b-1}}(x)}
\eea
for $b = 2, \dots , k$. From \Ref{FORMULA} it follows that the tuple
$\bs y^{k;\, c_1, \dots , c_{k}}$ corresponds to the rational section
\bea
S_{(c_1, \dots , c_k)} \ :\ x \ \mapsto\  f_{i_1, \dots , i_k} (g(x, c_1, \dots , c_k)) \times x
\eea
of the bundle $p'$. This section is horizontal with respect to the connection
$\nabla_{D_{\bs y^0}}'$, and we have
\bea
\phi_{\bs y^0,\, x} (\bs y^{k;\, c_1, \dots , c_{k}})\ =\
 f_{i_1, \dots , i_k} (g(x, c_1, \dots , c_k)) \ .
\eea
This means that
\bea
\phi_{\bs y^0, x_0} ( P^{[i_1, \dots , i_k]} )\ \subset \ B_{w^{-1}}\ .
\eea
It remains to show that every point in $B_{w^{-1}}$ is the limit of points
of $\phi_{\bs y^0, x_0} ( P^{[i_1, \dots , i_k]} )$, where the limit is taken in the sense of
the limit in Corollary \ref{Bruhat}, but that statement follows from

\begin{lemma} 
Let $x_0\in \C$ be such that $T_i(x_0)\neq 0$ and $y^0_i(x_0)\neq 0$ for 
$i = 1, \dots , r$. Then for any $(c_1^1, \dots , c_k^1) \in (\C - 0)^k$ there exists
a unique $(c_1^2, \dots , c_k^2) \in \C^k$ such that
\bea
(c_1^1, \dots , c_k^1) = g ( x_0; c_1^2, \dots , c_k^2) \ .
\eea
\hfill
$\square$
\end{lemma}

The proposition is proved.
\hfill
$\square$

Theorem \ref{second main} is a direct corollary of Proposition \ref{MAIN}.

\section{Solutions of Differential Equations}\label{solutions}

As we observed earlier, the Miura opers,
 associated with a population of critical points,
help to study the structure of the population.
In addition to that it
 turns out that for  a Miura oper $D$  associated with a critical point of a
population, 
all solutions of the differential equation $DY=0$ with values in the corresponding
group can be written explicitly in terms of critical points composing the population.
%In particular, every coordinate of every solution of the differential
%equation $D_{\bs y}\, Y\, =\, 0$ with values in a finite dimensional
%representation of $\tG$  can be written as a rational function
%of coordinates of tuples of the population of critical points originated at $\bs y$.

First we give formulas for solutions of the equation $DY=0$
for opers associated with Lie algebras of types $A_r,$ $ B_r$,
% $ C_N, D_N$ 
and then consider more general formulas for solutions 
which do not use the structure of the Lie algebra.

In this section $\g = \g(A)$ is a simple Lie algebra with Cartan matrix $A = (a_{i,j})$.

\subsection{Elimination of polynomials $T_i$}
Let $B = (b_{i,j})$ be the matrix inverse to $A$.

Let $D$ be the Miura $\tG$-oper   associated with
weights $\bs\La$, numbers $\bs z$, and a tuple $\bs y = (y_1, \dots , y_r)$.
Introduce polynomials $T_1(x), \dots , T_r(x)$ by formulas \Ref{T}.
Introduce a tuple of functions $\bs {\bar y} = (\bar y_1, \dots , \bar y_r)$ by
\bean\label{reduced polynomials}
\bar y_i\ =\ y_i \ \prod_{l=1}^r T_l^{-b_{i,l}} \ .
\eean

\begin{lemma}\label{reduced}
We have $D =  \partial  +  I + V$ where
\bean\label{V}
V \ = \ \sum_{j=1}^r  \ \mathrm{log}'(\bar y_j)\ H_j \ . 
\eean
\end{lemma}

\begin{proof} If $V$ is given by \Ref{V}, then
\bea
\langle \ {}^t\al_i, V \rangle  = 
\sum_{j=1}^r \mathrm{log}'(\bar y_j)\ \langle\ {}^t\al_i,\ {}^t\al^\vee_j  \rangle
=
  \sum_{j=1}^r \mathrm{log}'(\ \bar y_j^{\ {} a_{i,j}}\ )\ =\ 
\mathrm{log}'( \ T_i^{-1}\ \prod_{j=1}^r y_j^{\ {} a_{i,j}} \ ) \ .
\eea
\end{proof}

Assume that $\bs y$ represents a critical point of the master function
\Ref{master} associated with parameters $\bs z, \bs \La, \La_\infty$.

Let $\bs i = [i_1,  \dots , i_k], \ 
 i_j \in \{1, \dots , r\},$  
be a sequence of natural numbers. Let
$\bs y^{[i_1]} = ( y_1^{[i_1]}, \dots , y_r^{[i_1]})$,
$\bs y^{[i_1, i_2]} = ( y_1^{[i_1,i_2]}, \dots , y_r^{[i_1,i_2]})$,
\dots ,
$\bs y^{[i_1,  \dots , i_k]} = ( y_1^{[i_1, \dots , i_k]}, \dots , 
y_r^{[i_1, \dots , i_k]})$ be a sequence of tuples
associated with the critical point $\bs y$
and the sequence of indices $\bs i$, see Section \ref{general procedure}.
 Introduce functions
$\bar  y_{i}^{[i_1, \dots , i_l]}$ by 
\bean\label{reduced y}
\bar  y_{i}^{[i_1, \dots , i_l]} \ =\   y_{i}^{[i_1, \dots , i_l]}\ 
\prod_{l=1}^r T_l^{-b_{i,l}} \ .
\eean

\begin{lemma}
We have 
\bean\label{red rel}
W( \bar y_{i_1},  \bar y_{i_1}^{[i_1]})\ = 
 \ \prod_{j,\ j\neq i_1} \bar y_j^{- a_{i_1, j}}\ 
\eean
and $ \bar y_j^{[i_1]} = \bar y_j$ for $j\neq i_1$;
\ {} \  
for $l = 2, \dots , k$, we have
\bean\label{red rels}
% \label{ind eqn 2}
W( \bar y_{i_l}^{[i_1, \dots , i_{l-1}]},  \bar y_{i_l}^{[i_1, \dots , i_l]})\ = 
 \ \prod_{j,\ j\neq i_l}
(\bar y_j^{[i_1, \dots , i_{l-1}]})^{-a_{i_l, j}}\ 
\eean
and $ \bar y_{j}^{[i_1, \dots , i_{l}]} =  \bar y_{j}^{[i_1, \dots , i_{l-1}]}$
for $j\neq i_l$.
\hfill $\square$
\end{lemma}

The sequence of tuples
$\bs {\bar{y}}^{[i_1]} = ( \bar y_1^{[i_1]}, \dots ,  \bar y_r^{[i_1]})$,
$\bs  {\bar {y}}^{[i_1, i_2]} = ( \bar  y_1^{[i_1,i_2]}, \dots ,  \bar  y_r^{[i_1,i_2]})$,
\dots ,
\linebreak
$\bs { \bar {y}}^{[i_1, i_2, \dots , i_k]} = ( \bar  y_1^{[i_1,i_2, \dots , i_k]}, \dots , 
 \bar y_r^{[i_1,i_2, \dots , i_k]})$ will be called {\it the sequence of reduced tuples
associated with the critical point $\bs y$
and the sequence of indices $\bs i$}.

${}$

 The sequence of functions
$\bar y_{i_1}^{[i_1]}$, $\bar y_{i_2}^{[i_1,i_2]}$, \dots ,
$\bar y_{i_k}^{[i_1, i_2, \dots , i_k]}$, each defined up to multiplication by a non-zero
number, will be called {\it the reduced diagonal sequence of functions associated with
the critical point $\bs y$ and the sequence of indices $\bs i$.}

 The reduced diagonal sequence of functions
 determine the sequence of tuples $\bs {\bar y}^{[i_1]}$, 
$\bs {\bar y}^{[i_1,i_2]}$, \dots ,
$\bs {\bar y}^{[i_1, i_2, \dots , i_k]}$ uniquely.

%\bigskip

%Replacing the functions $y_i$ by functions $\bar y_i$ we achieve two goals: we express
%the oper $D$ in terms of the functions $\bar y_i$ only, and we express
%the relations of Lemma \ref{generation} in terms of 
%the functions $\bar y_i$ only.

%Our next goal is to solve the differential equation $DY = 0$ where $D$ 
%$D$ is the Miura $\tG$-oper   associated with weights $\bs\La = (\La_1, \dots , \Lambda_n)$, 
%numbers $\bs z = (z_1, \dots , z_n)$, and a tuple $\bs y = (y_1, \dots , y_r)$.

%\bigskip

In the next sections we  use the following  lemma.

\begin{lemma}\label{D0} 
${}$

\begin{enumerate}
\item[$\bullet$]
Consider the product $\prod_{j=1}^r \bar y_j^{\ - H_j}$ as a function of $x$ with values
in  the group $\tG$. Then
\bea
\prod_{j=1}^r \bar y_j^{\ - H_j}\ =\
\prod_{j=1}^r  \ y_j^{\ - H_j}
  T_j^{\ w_j}\ ,
\eea
where $w_1, \dots , w_r \in \ {}^t\h$ are the fundamental coweights, i.e.
$\langle \ {}^t\al_i , w_j \rangle = \delta_{i,j}$.
\item[$\bullet$]
 Let $D = \p + I + V$ be the Miura oper with $V$ given by formula \Ref{V}.
Define
\bea
\bar D \ =\ \p\ + \ \sum_{j=1}^r\ ( \ \prod_{l=1}^r\ \bar y_l^{\ -a_{j,l}}\ )\  F_j \ .
\eea
Then 
\bea D\ (\ \prod_{j=1}^r \bar y_j^{\ - H_j}) \ = \
(\ \prod_{j=1}^r \bar y_j^{\ - H_j}) \ \bar D \ .
\eea
\hfill $\square$
\end{enumerate}
\end{lemma}

\subsection{The $A_r$ critical points and $A_r$ opers}\label{AO} 
The Lie algebra $sl_{r+1}$ is of $A_r$-type. The Langlands dual to
$sl_{r+1}$ is $sl_{r+1}$. Let $F_1, \dots , F_r, H_1, \dots , H_r,
E_1, \dots , E_r$ be the Chevalley generators of $sl_{r+1}$. Let
$w_1, \dots , w_r$ be the fundamental coweights of $sl_{r+1}$.

We start with two examples.

Let $\g = sl_2$. Let $\bs y = (y_1)$ represent a critical point of the $sl_2$ master
function \Ref{master}  associated with parameters $\bs z, \bs \La, \La_\infty$.
Introduce the function $\bar y_1 = y_1 T_1^{\ -1/2}$, see \Ref{reduced polynomials}.
Then the Miura oper associated with $\bs y$ has the form
\bea
D = \p + F_1 + \mathrm{log}'(\bar y_1) H_1\ .
\eea
Let $\bar y_1^{[1]}$ be the reduced diagonal sequence
of functions associated with $\bs y$ and the sequence of indices $[1]$, \ 
in other words,
$W(\bar y_1, \bar y_1^{[1]}) = 1$. Then
\bea
Y \ =\ \bar y_1^{\ - H_1}\ e^ { \frac{\bar y_1^{[1]}}{\bar y_1} F_1}\ 
\eea
is a solution of the differential equation $DY=0$ with values in $\mathrm{SL}\,(2,\C)$.
Indeed,
\bea
DY = \bar y_1^{\ - H_1}  (\p + \frac 1{(\bar y_1)^2} F_1 )\
e^ { \frac{\bar y_1^{[1]}}{\bar y_1} F_1}\ =
Y (\p + (\left( \frac{\bar y_1^{[1]}}{\bar y_1}\right)' 
\!\!+ \frac 1{(\bar y_1)^2}) F_1 )\ \mathrm{id} 
=\ Y\ \p \ \mathrm{id} = 0  .
\eea

Let $\g = sl_3$. Let $\bs y = (y_1,y_2)$ represent a critical point of the
$sl_3$ master
function \Ref{master}  associated with parameters $\bs z, \bs \La, \La_\infty$.
Introduce the functions $\bar y_1 = y_1 T_1^{\ -2/3}T_2^{\ -1/3}$, 
$\bar y_2 = y_2 T_1^{\ -1/3}T_2^{\ -2/3}$, see \Ref{reduced polynomials}.
Then the Miura oper associated with $\bs y$ has the form
\bea
D = \p + F_1 + F_2 + \mathrm{log}'(\bar y_1) H_1 + \mathrm{log}'(\bar y_2) H_2\ .
\eea
Let $\bar y_1^{[1]}, \bar y_2^{[1,2]}$ be the reduced diagonal sequence
of functions associated with $\bs y$  and the sequence of indices $[1,2]$, \ 
in other words,
\bea
W(\bar y_1, \bar y_1^{[1]}) = \bar y_2 , 
\qquad
W(\bar y_2, \bar y_2^{[1,2]}) = \bar y_1^{[1]} .
\eea
Let $\bar y_2^{[2]}$ be the reduced diagonal sequence
of functions associated with $\bs y$  and the sequence of indices $[2]$, \
in other words,
$W(\bar y_2, \bar y_2^{[2]}) = \bar y_1$.  Then
\bea
Y\ =\ \bar y_1^{\ - H_1}\ \bar y_2^{\ - H_2}\ e^{ \frac{\bar y_1^{[1]}}{\bar y_1} F_1}
\ e^{ \frac{\bar y_2^{[1,2]}}{\bar y_2} [F_2,F_1]}\
e^{ \frac{\bar y_2^{[2]}}{\bar y_2} F_2}
\eea
is a solution of the differential equation $DY=0$ with values in $\mathrm{SL}\,(3,\C)$.
Indeed, by Lemma \ref{D0} it suffices to show that
\bea
\bar Y\ = \
\ e^{ \frac{\bar y_1^{[1]}}{\bar y_1} F_1}
\ e^{ \frac{\bar y_2^{[1,2]}}{\bar y_2} [F_2,F_1]}
\ e^{ \frac{\bar y_2^{[2]}}{\bar y_2} F_2}
\eea
is a solution of the differential equation $\bar D Y=0$ where
\bea
\bar D \ =\ \p + \frac {\bar y_2}{(\bar y_1)^2} F_1 + \frac {\bar y_1}{(\bar y_2)^2} F_2 \ .
\eea
Indeed, 
\bea
\bar D \bar Y  & = &  e^{ \frac{\bar y_1^{[1]}}{\bar y_1} F_1}
(\p + (\left( \frac{\bar y_1^{[1]}}{\bar y_1}\right)' \!\!
+ \frac 1{(\bar y_1)^2}) F_1  + 
\frac {\bar y_1^{[1]}}{(\bar y_2)^2} [F_2,F_1] + \frac {\bar y_1}{(\bar y_2)^2} F_2 )
\ e^{ \frac{\bar y_2^{[1,2]}}{\bar y_2} [F_2,F_1]}
\ e^{ \frac{\bar y_2^{[2]}}{\bar y_2} F_2}
\\
& = & 
e^{ \frac{\bar y_1^{[1]}}{\bar y_1} F_1}
(\p + 
\frac {\bar y_1^{[1]}}{(\bar y_2)^2} [F_2,F_1] + \frac {\bar y_1}{(\bar y_2)^2} F_2 )
\ e^{ \frac{\bar y_2^{[1,2]}}{\bar y_2} [F_2,F_1]}
\ e^{ \frac{\bar y_2^{[2]}}{\bar y_2} F_2}
\\
&=&
e^{ \frac{\bar y_1^{[1]}}{\bar y_1} F_1} e^{ \frac{\bar y_2^{[1,2]}}{\bar y_2} [F_2,F_1]}
(\p + 
(\left( \frac{\bar y_2^{[1,2]}}{\bar y_2}\right)' +
\frac {\bar y_1^{[1]}}{(\bar y_2)^2} ) [F_2,F_1] + \frac {\bar y_1}{(\bar y_2)^2} F_2 )
\ e^{ \frac{\bar y_2^{[2]}}{\bar y_2} F_2}
\eea
\bea
&=&
e^{ \frac{\bar y_1^{[1]}}{\bar y_1} F_1} e^{ \frac{\bar y_2^{[1,2]}}{\bar y_2} [F_2,F_1]}
(\p + 
 \frac {\bar y_1}{(\bar y_2)^2} F_2 )
\ e^{ \frac{\bar y_2^{[2]}}{\bar y_2} F_2}
\\
&=& 
e^{ \frac{\bar y_1^{[1]}}{\bar y_1} F_1} e^{ \frac{\bar y_2^{[1,2]}}{\bar y_2} [F_2,F_1]}
e^{ \frac{\bar y_2^{[2]}}{\bar y_2} F_2}
(\p + 
 (\left( \frac{\bar y_2^{[2]}}{\bar y_2}\right)' +
\frac {\bar y_1}{(\bar y_2)^2} ) F_2 )\ {} \mathrm{id} 
\\
&=&
e^{ \frac{\bar y_1^{[1]}}{\bar y_1} F_1} e^{ \frac{\bar y_2^{[1,2]}}{\bar y_2} [F_2,F_1]}
e^{ \frac{\bar y_2^{[2]}}{\bar y_2} F_2}
\ \p  \ \mathrm{id } = 0 \ .
 \eea

Now consider the general case.
Let $\g = sl_{r+1}$.
Let $\bs y = (y_1, \dots , y_r)$ represent a critical point of the $sl_{r+1}$
master
function \Ref{master}  associated with parameters $\bs z, \bs \La, \La_\infty$.
Introduce the functions $\bar y_1, \dots , \bar y_r$ by formula  
\Ref{reduced polynomials}, where $B=(b_{i,j})$ is the matrix inverse to the
Cartan matrix of $sl_{r+1}$.
Then the Miura oper associated with $\bs y$ has the form
\bean\label{oper SL}
D\ =\ \p\ +\ \sum_{j=1}^r\ F_j\ +\ \sum_{j=1}^r \ \mathrm{log}'(\bar y_j)\  H_j\ .
\eean
For $i = 1, \dots , r$,\ let 
$ y_i^{[i]},\  y_{i+1}^{[i,i+1]},\ \dots ,\
 y_{r}^{[i,\dots , r]}$  be the  diagonal sequence
of polynomials associated with $\bs y$  and the sequence of indices $[i,i+1, \dots , r]$, \ 
in other words,
\bea
W( y_i,  y_i^{[i]}) \ =\ T_i\ y_{i-1} \ y_{i+1}  ,  
\qquad
W( y_{i+1},  y_{i+1}^{[i,i+1]}) \ =\ T_{i+1}\ y_i^{[i]} \bar y_{i+2} ,
\ {}\ \dots ,
\eea
\bea
W( y_{r-1},  y_{r-1}^{[i,\dots , r-1]}) \ =\ 
T_{r-1}\ y_{r-2}^{[i, \dots , r-2]}\  y_{r} ,
\qquad
W( y_{r},  y_{r}^{[i,\dots , r]}) \ =\ 
T_r\ y_{r-1}^{[i, \dots , r-1]} .
\eea
Define $r+1$ functions $Y_0, Y_1, \dots , Y_r$
of $x$ with values in $\mathrm{SL}\,(r+1,\C)$ by the formulas
\bea
Y_0\ = \ \prod_{j=1}^r\  y_j ^{\ -H_j} T_j^{\ w_j} , 
\qquad
Y_i\ = \ \prod_{j=i}^r \ e^{ \frac { y_j^{[i,\dots ,j]}}{ y_j}
\ [F_j,[F_{j-1},[...,[F_{i+1},F_i]...]]]} ,
\qquad 
\text{for}\ i > 0 .
\eea
Note that inside each product the factors commute.

\begin{theorem}\label{sl}
The product $Y_0 Y_1 \dots Y_r$ is a solution of the differential equation
$DY = 0$ with values in $\mathrm{SL}\,(r+1,\C)$ where $D$ is given by \Ref{oper SL}.
\end{theorem}

Note that if $Y(x)$ is a solution of the equation $DY = 0$ and $g\in \mathrm{SL}\,(r+1, \C)$,
then $Y(x)g$ is a solution too.

The proof of the theorem is straightforward. One uses Lemma \ref{D0} and then shows that
\bea
( \p\ + \ \sum_{j=i}^r\ \frac{\bar y_{j-1} \bar y_{j+1}}{\bar y_j^{\ 2}}\ F_j )
\ Y_i\ = \ Y_i \
( \p\ + \ \sum_{j=i+1}^r\ \frac{\bar y_{j-1} \bar y_{j+1}}{\bar y_j^{\ 2}}\ F_j )
\eea
for $i = 1 , \dots , r$. In this formula we set $\bar y_0 = \bar y_{r+1} = 1$.

\subsection{The $B_r$ critical points and $C_r$ opers}\label{BC} 
Consider the root system of type $B_r$. Let
$\al_1, \dots , \al_{r-1}$ be the long simple roots and $\al_r$ the short one.
 We have 
\bea
(\al_r , \al_r)\ =\ 2 , \qquad
(\al_i , \al_i)\ =\ 4 , \qquad
(\al_i , \al_{i+1}) \ =\ - 2 , \qquad i = 1, \dots , r-1 ,
\eea
and all other scalar products are equal to zero.
The root system $B_r$ corresponds to the Lie algebra $so_{2r+1}$.
Let $\h_B$ be its Cartan subalgebra. 
%The spinor group $\mathrm{Spin}(2r+1,\C)$ is the simply connected group
%with Lie algebra $so_{2r+1}$.

\bigskip

Consider the root system of type $C_r$. 
The root system $C_r$ corresponds to the Lie algebra $sp_{2r}$.
Let $F_1, \dots , F_r,\ H_1, \dots , H_r,\ E_1, \dots , E_r$ \ be 
its Chevalley generators and $w_1, \dots , w_r$ the fundamental coweights.
The symplectic group $\mathrm{Sp}(2r,\C)$ 
is the simply connected group with Lie algebra $sp_{2r}$.

The Lie algebras $so_{2r+1}$ and $sp_{2r}$ are Langlands dual.

\bigskip

We consider also the root system of type $A_{2r-1}$ with simple roots $\al^A_1, \dots ,
\al^A_{2r-1}$. The root system $A_{2r-1}$ corresponds to the Lie algebra $sl_{2r}$. We
denote $\h_A$ its Cartan subalgebra.

We have a map $\h^*_B \to \h_A^*, \  \La \mapsto \La^A$, where $\La^A$ is defined by
\be
\langle \La^A\  , \ (\al_i^A)^\vee \rangle\ =\ \langle \La^A\ ,\ (\al_{2r-i}^A)^\vee
\rangle \ =\
\langle \La \ , \ (\al_i)^\vee \rangle , \qquad i = 1 , \dots ,  r . 
\ee

Let $\La_1, \dots, \La_n \in \h^*_B$ be dominant integral $so_{2r+1}$-weights,
$z_1, 
\dots , z_n$ complex numbers. Let the polynomials $T_1, \dots , T_r$ be given by \Ref{T}.
Remind that an $r$-tuple
  of polynomials $\bs y$ represents a critical point of a master function
associated with $so_{2r+1}$, \ $\La_1, \dots , \La_n$, \ $z_1, \dots , z_n$,
if and only if $\bs y$ is generic with respect to weights
\ $\La_1, \dots , \La_n$ of $so_{2r+1}$, and points $z_1, \dots z_n$
and
there exist polynomials $\tilde y_i$,\  $i = 1, \dots , r$,  such that
\bea
&& 
W(y_i \ ,\  \tilde y_i )\ =\ T_i\ y_{i-1}\ y_{i+1} ,\qquad i\ =\ 1, \dots , r-1 ,
\\
&&
W(y_r \ , \ \tilde y_r)\ =\ T_r\ y_{r-1}^2 .
\eea

For an $r$-tuple of polynomials $\bs y = (y_1, \dots, y_r)$, let 
$\bs u$ be the $2r-1$-tuple of polynomials
$(u_1, \dots , u_{2r-1}) = (y_1, \dots , y_{r-1}, y_r, y_{r-1}, \dots, y_1)$.

\begin{lemma}[\cite{MV1}]
An $r$-tuple $\bs y$ represents a critical point of the
 $so_{2r+1}$  master function 
associated with  $\La_1, \dots , \La_n$, \ $z_1, \dots , z_n$,
if and only if the $2r-1$-tuple of polynomials
$\bs u$ represents a critical point of the $sl_{2r}$ master function
associated with \ $\La_1^A, \dots ,$ $ \La_n^A$, \linebreak
 $z_1, $ $\dots ,$ $ z_n$.
\hfill $\square$
\end{lemma}

\bigskip

We start with an example. 
Let $\bs y = (y_1, y_2)$ represent a critical point of the $so_{3}$
master function \Ref{master}  associated with parameters $\bs z, \bs \La, \La_\infty$.
Set $\bar y_1 = y_1 T_1^{\ -1} T_2^{\ -1/2}, \linebreak
 \bar y_2 = y_2 
T_1^{\ -1} T_2^{\ -1}$, see \Ref{reduced polynomials}. 

Let $\bs u = (u_1, u_2, u_3) = (y_1, y_2, y_1)$ be the tuple representing the corresponding
$sl_4$ critical point. Set $\bs {\bar u} = (\bar u_1, \bar u_2, \bar u_3) 
= (\bar y_1, \bar y_2, \bar y_1)$.

Let $\bar y_1^{[1]}, \bar y_2^{[1,2]}$ be the $so_{3}$
reduced diagonal sequence
of functions associated with $\bs y$  and the sequence of indices $[1,2]$, \
in other words,
\bea
W(\bar y_1, \bar y_1^{[1]}) = \bar y_2\ ,
\qquad  
W(\bar y_2, \bar y_2^{[1,2]}) = (\bar y_1^{[1]})^{ 2}\ .
\eea
Let $\bar y_2^{[2]}$ be the $so_{3}$ reduced diagonal sequence
of functions associated with $\bs y$  and the sequence of indices $[2]$, \ 
in other words, 
\bea
W(\bar y_2, \bar y_2^{[2]}) = (\bar y_1)^{2}\ .
\eea
Let $\bar u_1^{[1]} = \bar y_1^{[1]},\  \bar u_2^{[1,2]}$ 
be the $sl_4$ reduced diagonal sequence
of functions associated with $\bs u$  and the sequence of indices $[1,2]$, \
in other words,
\bea
W(\bar y_1, \bar y_1^{[1]}) = \bar y_2\ ,
\qquad  
W(\bar y_2, \bar u_2^{[1,2]}) = \bar y_1^{[1]}\bar y_1 \ .
\eea

Then
\bea
Y\ =\ \bar y_1^{\ - H_1}\ \bar y_2^{\ - H_2}\ e^{ \frac{\bar y_1^{[1]}}{\bar y_1} F_1}
e^{ \frac 12\frac{\bar y_2^{[1,2]}}{\bar y_2} [[F_2,F_1],F_1]}\
 e^{ \frac{\bar u_2^{[1,2]}}{\bar y_2} [F_2,F_1]}\
e^{ \frac{\bar y_2^{[2]}}{\bar y_2} F_2}
\eea
is an $\mathrm{Sp}\,(4,\C)$-valued solution of the 
differential equation $DY = 0$ where
\bea
D\ =\  \p\ +\ F_1\ +\ F_2\ +\
\mathrm{log}'(\bar y_1)\, H_1\ +\ \mathrm{log}'(\bar y_2) \,H_2\ .
\eea 

Indeed, denote the factors of $Y$ by $P_1, \dots , P_6$ counting from the left.
By Lemma \ref{D0} it suffices to show that the product $P_3P_4P_5P_6$ is a solution
of the equation $\bar D Y = 0$ where
\bea
\bar D = \p + 
\frac{\bar y_2}{\bar y_2 ^{\ 2}} F_1 +
\frac{\bar y_1^{\ 2}}{\bar y_2 ^{\ 2}} F_2 \ .
\eea
We have
\bea
 \bar D \ P_3P_4P_5P_6  
&=&
P_3\ (\p + 
\frac{\bar y_1 \bar y_1^{[1]}}{\bar y_2 ^{\ 2}} [F_2,F_1] +
\frac 12 \frac{ (\bar y_1^{[1]})^{ 2}}{\bar y_2 ^{\ 2}} [[F_2,F_1],F_1] +
\frac{\bar y_1^{\ 2}}{\bar y_2 ^{\ 2}} F_2 )\ P_4P_5P_6 
\\
&=&
P_3P_4\ (\p + 
\frac{\bar y_1 \bar y_1^{[1]}}{\bar y_2 ^{\ 2}} [F_2,F_1] +
%\frac 12 \frac{ (\bar y_1^{[1]})^{ 2}}{\bar y_2 ^{\ 2}} [[F_2,F_1],F_1] +
\frac{\bar y_1^{\ 2}}{\bar y_2 ^{\ 2}} F_2 )\ P_5P_6 
\\
&=&
P_3P_4P_5\ (\p + 
%\frac{\bar y_1 \bar y_1^{[1]}}{\bar y_2 ^{\ 2}} [F_2,F_1] +
%\frac 12 \frac{ (\bar y_1^{[1]})^{ 2}}{\bar y_2 ^{\ 2}} [[F_2,F_1],F_1] +
\frac{\bar y_1^{\ 2}}{\bar y_2 ^{\ 2}} F_2 )\ P_6 
\ =\
P_3P_4P_5P_6\ \p \ \mathrm{id}\ =\ 0\ .
\eea

\bigskip

Now consider the general case.
Let $\bs y = (y_1, \dots , y_r)$ represent a critical point of the $so_{2r+1}$
master function \Ref{master}  associated with parameters $\bs z, \bs \La, \La_\infty$.
Introduce the functions $\bar y_1, \dots , \bar y_r$ by formula  
\Ref{reduced polynomials}, where $B=(b_{i,j})$ is the matrix inverse to the
Cartan matrix of $so_{2r+1}$.
Then the $sp_{2r}$ Miura oper associated with $\bs y$ has the form
\bean\label{so 2r + 1}
D\ =\ \p\ +\ \sum_{j=1}^r\ F_j \ +\ \sum_{j=1}^r\ \mathrm{log}'(\bar y_j)\  H_j\ .
\eean

Let $\bs u = (u_1, \dots , u_{2r-1}) = (y_1,\dots, y_r,\dots, y_1)$ 
be the tuple representing the corresponding
$sl_{2r}$ critical point.
% Set $\bs {\bar u} = (\bar u_1, \dots,  \bar u_{2r-1}) 
%= (\bar y_1, \dots , \bar y_r, \dots , \bar y_1)$.

For $i = 1, \dots , r$,\ let 
$ y_i^{[i]}$,\ $ y_{i+1}^{[i,i+1]}$,\ \dots ,\
$ y_{r}^{[i,\dots , r]}$  be the $so_{2r+1}$ diagonal 
%\linebreak 
sequence
of polynomials
 associated with $\bs y$  and the sequence of indices $[i,i+1, \dots , r]$, \ 
in other words,
\bea
W( y_i,  y_i^{[i]}) \ =\ T_i\ y_{i-1} \  y_{i+1}  ,  
\qquad
W( y_{i+1},  y_{i+1}^{[i,i+1]}) \ =\ T_{i+1}\ y_i^{[i]}\  y_{i+2} ,
\ {} \dots ,
\eea
\bea
W( y_{r-1},  y_{r-1}^{[i,\dots , r-1]}) \ =\ 
T_{r-1}\ y_{r-2}^{[i, \dots , r-2]}\  y_{r} ,
\qquad
W( y_{r},  y_{r}^{[i,\dots , r]}) \ =\ 
T_r\ ( y_{r-1}^{[i, \dots , r-1]})^{2} .
\eea
For $i = 1, \dots , r-1$,\ let 
\bea
 u_i^{[i]} =  y_i^{[i]}, 
\qquad
 u_{i+1}^{[i,i+1]} =  y_{i+1}^{[i,i+1]}, 
\qquad
 \dots\ {} ,
\qquad
 u_{r-1}^{[i,\dots , r-1]} =  y_{r-1}^{[i,\dots , r-1]},
\eea
\bea
 u_{r}^{[i,\dots , r]} ,
\qquad
 u_{r+1}^{[i,\dots , r+1]} , 
\qquad
 \dots \  ,
\qquad
 u_{2r-i-1}^{[i, \dots , 2r-i-1]}
\eea
 be the $sl_{2r}$ diagonal sequence
of polynomials associated with $\bs u$  and the sequence of indices 
$[i,i+1, \dots , 2r-i-1]$, \ 
in other words,
\bea
W( y_{r},  u_{r}^{[i,\dots , r]}) \ =\ 
T_r\ y_{r-1}^{[i, \dots , r-1]}\  y_{r-1} ,
\qquad
W( y_{r-1},  u_{r+1}^{[i,\dots , r+1]}) \ =\ 
T_{r-1}\ u_{r}^{[i, \dots , r]} \  y_{r-2} \ ,
\eea
\bea
W( y_{r-l},  u_{r+l}^{[i,\dots , r+l]}) \ =\ 
T_{r-l}\ u_{r+l-1}^{[i, \dots , r+l-1]}\  y_{r-l-1} \ ,
\qquad 
\text{for} \ {}\ 
l = 2 , \dots , r-i-1 .
\eea

For $i \in \{1, \dots , r\}$ set $F_{i,i} = F_i$. 
For $1\leq i < j < r$ set
\bea
F_{i,j}\ = \ [F_j,[F_{j-1},[...,[F_{i+1},F_i]...]]]\ .
\eea
 Set $F^*_{i,r} = [F_r, F_{i,r-1}]$ and for
$1\leq i < j < r$ set
\bea
F_{i,j}^*\ =\ [F_j,[F_{j+1},[...[F_{r-2},[F_{r-1},F_{i,r}^*]]...]]]\ .
\eea

Define $r+1$ functions $Y_0, Y_1, \dots , Y_r$
of $x$ with values in $\mathrm{Sp}\,(2r,\C)$ by the formulas
\bea
Y_0 \ =\  \prod_{j=1}^r\  y_j ^{\ -H_j} T_j^{\ w_j}\ ,
\eea
\bea
Y_i\  =\  \left( \ \prod_{j=i}^{r-1} \ e^{ \frac { y_j^{[i,\dots ,j]}}{ y_j}
\ F_{i,j}} \ \right)\
e^{ \frac 12 \frac{ y_r^{[i, \dots , r]}}{ y_r}\ [[F_r,F_{i,r-1}],F_{i,r-1}]}\
\left(\ \prod_{j=r}^{2r-i-1}
e^{ \frac { u_j^{[i,\dots ,j]}}{ y_{2r - j}}\ F^*_{i,2r-j}} \ \right)
\eea
for $i \in \{1, \dots, r-1\}$, and\  $Y_r = 
e^{ \frac { y_r^{[r]}}{ y_{r}}\ F_r}$.

Note that inside each product the factors commute.

\begin{theorem}\label{bc}
The product $Y_0 Y_1 \dots Y_r$ is a solution of the differential equation
$DY = 0$ with values in $\mathrm{Sp}\,(2r,\C)$ where $D$ is given by \Ref{so 2r + 1}.
\end{theorem}

The proof is straightforward. One uses Lemma \ref{D0} and then shows that
\bea
( \p\ + \ \frac {\bar y_{r-1}^{\ 2}}{\bar y_{r}^{\ 2}} F_r
\ +\ \sum_{j=i}^{r-1}\ \frac{\bar y_{j-1} \bar y_{j+1}}{\bar y_j^{\ 2}}\ F_j )
\ Y_i\ = \ Y_i \
( \p\  + \ \frac {\bar y_{r-1}^{\ 2}}{\bar y_{r}^{\ 2}} F_r
\ +\ \sum_{j=i+1}^{r-1}\ \frac{\bar y_{j-1} \bar y_{j+1}}{\bar y_j^{\ 2}}\ F_j )
\eea
for $i = 1 , \dots , r-1$. 
%In this formula we set $\bar y_0 = \bar y_{r+1} = 1$.

\bigskip

{\bf Remark.} Theorems \ref{sl} and \ref{bc} give explicit formulas for
solutions of the differential equation $DY = 0$ where $D$ is the Miura oper
associated to a critical point of type $A_r$ or $B_r$. In a similar way one can construct
explicit formulas for solutions in the case of the Miura oper associated to a critical point
of type $C_r$, cf. Section 7 in \cite{MV1}.

\subsection{General formulas for solutions} 
Let $\g$ be a simple Lie algebra with Cartan matrix $A$. Let
${}^t\g$ be its Langlands dual with Chevalley generators 
 $F_1, \dots , F_r,\ H_1, \dots , H_r$, \linebreak
$ E_1, \dots , E_r$. Let $w_1, \dots , w_r$ be the fundamental coweights of ${}^t\g$.
 Let $\tG$ be the complex simply connected
 Lie group with Lie algebra $\tg$.

Let $V$ be a complex finite dimensional representation of $\tG$. Let $v_{low}$ be 
a lowest weight
vector of $V$, \ {} ${}^t\n_-\, v_{\mathrm{low}} = 0$.

Let $\bs y = (y_1, \dots , y_r)$ represent a critical point of the $\g$
master function \Ref{master}  associated with parameters $\bs z, \bs \La, \La_\infty$.
Let $D_{\bs y}$ be the ${}^t\g$ Miura oper associated with $\bs y$.

%Introduce the functions $\bar y_1, \dots , \bar y_r$ by formula  
%\Ref{reduced polynomials}, where $B=(b_{i,j})$ is the matrix inverse to the
%Cartan matrix of $\g$.
%Then the ${}^t\g$ Miura oper associated with $\bs y$ has the form
%\bea
%D_{\bs y}\ =\ \p\ +\ \sum_{j=1}^r\ F_j \ +\ \sum_{j=1}^r\ \mathrm{log}'(\bar y_j)\  H_j\ .
%\eea

Let $\bs i = [i_1,  \dots , i_k], \ 
 i_j \in \{1, \dots , r\},$   
be a sequence of natural numbers. Let
$\bs y^{[i_1]} = ( y_1^{[i_1]}, \dots , y_r^{[i_1]})$,
$\bs y^{[i_1, i_2]} = ( y_1^{[i_1,i_2]}, \dots , y_r^{[i_1,i_2]})$,
\dots ,
$\bs y^{[i_1,  \dots , i_k]} = ( y_1^{[i_1, \dots , i_k]}, \dots , 
y_r^{[i_1, \dots , i_k]})$ be a sequence of tuples
associated with the critical point $\bs y$
and the sequence of indices $\bs i$, see Section \ref{general procedure}.
%Introduce functions $\bar  y_1^{[i_1, \dots , i_k]}, \dots , 
%\bar y_r^{[i_1, \dots , i_k]}$ by formulas \Ref{reduced polynomials}.

\begin{theorem}\label{Solutions}
The $V$-valued function
\bea\label{Sol}
Y\ = \ \exp\!\left( -  \mathrm{log}'  \left(
\frac{ y_{i_1}^{[i_1]}}{ y_{i_1}}\right) E_{i_1} \right)
\exp\!\left( -  \mathrm{log}'  \left( \frac{ y_{i_2}^{[i_1,i_2]}}{ y_{i_2}^{[i_1]}}
\right) E_{i_2} \right)
\ \dots \phantom{aaaaaaaaaaaaaaaaaaaaaa}
\eea
\bea
\phantom{aaaaaaaaaaaaa}
 \exp\!\left(  - \mathrm{log}'  \left( 
\frac{ y_{i_k}^{[i_1,\dots , i_k]}}{ y_{i_k}^{[i_1, \dots , i_{k-1}]}}\right)
 E_{i_k} \right)\ {}
 \prod_{j=1}^r \ 
( \,y_{j}^{[i_1, \dots , i_{k}]}\,)^{- H_j} T_j^{\ w_j} \ v_{\mathrm{low}}
\eea
is a solution of the differential equation 
\bea
D_{\bs y}\ Y\ =\ 0\ .
\eea
\end{theorem}

The proof is straightforward and follows from the identity
\bea
 D_{\bs y^{[i_1,  \dots , i_j]}} \ =\
 \exp\!\left( \mathrm{ad}\,\mathrm{log}'  \left( 
\frac{ y_{i_j}^{[i_1,\dots , i_j]} } { y_{i_j}^{[i_1, \dots , i_{j-1}]} }
\right) E_{i_j} \right) \cdot D_{ \bs y^{[i_1, \dots , i_{j-1}]} }\ , 
% \exp\!\left(  - \mathrm{log}'  \left( 
%\frac{ y_{i_j}^{[i_1, \dots , i_j]} }{ y_{i_j}^{[i_1, \dots , i_{j-1}]}}\right) E_{i_j} 
%\right)
% = 
% D_{\bs y^{[i_1,  \dots , i_j]}} \,  ,
\eea 
see Theorem \ref{main}.

Let $d$ be the determinant of the Cartan matrix of $\g$.

\begin{corollary}\label{rationality}
Every coordinate of every solution of the equation $D_{\bs y}\, Y\,
=\, 0$ with values in a finite dimensional representation of \ $\tG$
can be written as a rational function $R(f_1, \dots , f_N; \
T_1^{1/d}, \dots , T_r^{1/d}\,)$ of functions $ T_1^{1/d}, \dots ,
T_r^{1/d}$ and suitable polynomials $f_1, \dots , f_N$ which appear as
coordinates of tuples in the $\g$ population $P_{\bs y}$ originated at
$\bs y$.
\end{corollary}

Since $\tG$ has a faithful finite dimensional representation, the
solutions of the differential equation $D_{\bs y}\, Y\, =\, 0$ with
values in $\tG$ also can be written as rational functions of functions
$ T_1^{1/d}, \dots , T_r^{1/d}$ and coordinates of tuples of $P_{\bs
y}$, cf. Sections \ref{AO} and \ref{BC}.

\end{document}